\documentclass[11pt,leqno]{article}
\usepackage{amsthm,amsfonts,amssymb,epsfig,graphics,amsmath,oldgerm}
   \usepackage[latin1]{inputenc}\relax

\newcommand{\RR}{{\mathbb R}}
\newcommand{\NN}{{\mathbb N}}
\newcommand{\ZZ}{{\mathbb Z}}

\newcommand{\cC}{\mathcal{C}}
\newcommand{\cH}{\mathcal{H}}
\newcommand{\cL}{\mathcal{L}}

\newcommand{\eps}{\varepsilon }
\newcommand{\vp}{\varphi}
\newcommand{\D}{\partial }

\newcommand{\na}{{\nabla}}
\newcommand{\Log}{\mathrm{Log}}
\newcommand{\la}{\langle}
\newcommand{\ra}{\rangle}
\newcommand{\mez}{\frac{1}{2}}
\newcommand{\ml}{\frac{1}{2} log}
\newcommand{\re}{{\mathrm{Re }}}

\newtheorem{theo}{Theorem}[section]
\newtheorem{prop}[theo]{Proposition}
\newtheorem{cor}[theo]{Corollary}
\newtheorem{lem}[theo]{Lemma}
\newtheorem{defi}[theo]{Definition}
\newtheorem{ass}[theo]{Assumption}

\newtheorem{rem}[theo]{Remark}

\numberwithin{equation}{section}
\title{The Cauchy  Problem for Wave Equations with non Lipschitz Coefficients}

   \author{Ferruccio Colombini \footnote{Dipartimento di Matematica,
Universit\`a di Pisa, Largo B. Pontecorvo 5, 56127 Pisa, Italia.
Email : Colombini@dm.unipi.it }
   and Guy M\'etivier \footnote{MAB,  Universit\'e de Bordeaux I,
33405 Talence cedex, France. Email : Guy.Metivier@math.u-bordeaux.fr} }

\newcount\madechanges
\def\caution#1{\ifnum \madechanges=1 \affixmessage{#1}%
\else \relax \fi}
\def\affixmessage#1{\marginpar{{\footnotesize \em #1} \openup
     -.3\baselineskip }}
\begin{document}

\maketitle

\begin{abstract}
In this paper we study the Cauchy problem for second order strictly 
hyperbolic operators of the form
\begin{equation*}
 L u := \sum_{j, k 
= 0}^n \D_{y_j} \big( a_{j, k} \D_{y_k} u \big)  +
   \sum_{j=0}^n \{   b_j \D_{y_j} u + \D_{y_j} (  c_j u)\}    + d u  = f,
\end{equation*}
 when the coefficients of the principal 
part are not Lipschitz continuous, but only
``Log-Lipschitz'' with 
respect to all the variables. This class of equation is invariant under changes of variables
and therefore suitable for a local analysis. 
In particular, we show local existence, local uniqueness and finite speed of propagation for the  noncharacteristic Cauchy problem. 
\end{abstract}
\tableofcontents

\section{Introduction}

In this paper we study the well-posedness of the Cauchy problem for
second order
strictly hyperbolic equations whose coefficients are not Lipschitz continuous:
\begin{equation}
\label{maineq}
   L u := \sum_{j, k = 0}^n \D_{y_j} \big( a_{j, k} \D_{y_k} u \big)  +
   \sum_{j=0}^n \{   b_j \D_{y_j} u + \D_{y_j} (  c_j u)\}    + d u  = f .
\end{equation}
This question has already been studied
in the case that the second order  part has the special form,  in coordinates
$ y = (t, x)$:
\begin{equation}
\label{speceq}
\D_t^2   -   \sum_{j, k = 1}^n \D_{x_j} \big( a_{j, k} \D_{x_k} u \big)
\end{equation}
and the Cauchy data are given on the space-like  hyperplane $\{ t = 0 \}$.
In this case,  when the  coefficients  depend only
on the time variable $t$, F. Colombini, E. De Giorgi and  S. Spagnolo
(\cite{CGS})
have proved that  the Cauchy problem
   is in general ill-posed in $C^\infty$  when the coefficients are
only H\"older continuous of order $\alpha < 1$,  but is  well-posed
in appropriate Gevrey spaces.
This has been extended to the case where the coefficients are
H\"older in time and Gevrey in $x$
(\cite{Ni, Ja}).
Moreover, it is also proved in  \cite{CGS} that the Cauchy problem is
well posed in
$C^\infty$ when the coefficients,  which depend only on time,
are  ``Log-Lischitz'' (in short  LL) : recall that a function $a$
of variables $y$  is said to be LL on a domain
$\Omega$ if there is a constant $C$ such that
\begin{equation}
\label{LL}
| a(y) - a(y') | \le C | y - y' | \, \Big( 1 +  \big| \mathrm{Log}
|y - y' | \big| \Big)
\end{equation}
for all $y$ and $y'$ in $\Omega$.
In \cite{CGS}, it is proved that for LL coefficients depending only
on $t $  and for initial data in the Sobolev spaces $H^s \times
H^{s-1}$, the solution satisfies
\begin{equation}
\label{loss}
u(t, \cdot ) \in H^{s- \lambda t} ,   \quad  \D_t u(t, \cdot ) \in
H^{s-1- \lambda t}
\end{equation}
with $\lambda $ depending only on the $LL$ norms of the coefficients
and the constants of hyperbolicity.
In particular, there is a loss of smoothness as   time evolves and
    this loss does occur in general when the coefficients are not
Lipschitz continuous,
    and is sharp, as shown in \cite{CC}.

The analysis of the  $C^\infty$ well-posedness has been extended
by F. Colombini and N. Lerner (\cite{CL}) to the case of equations,
still with principal
part \eqref{speceq},  whose  coefficients also depend  on the space
variables $x$.
   They show that the Cauchy problem is well-posed if the coefficients
are LL in time and
   $C^\infty$ in $x$.  They also study the problem under the natural
assumption of
   isotropic LL smoothness in $(t, x)$.  In this case one has to
multiply  LL functions with
   distributions in $H^s$. This is
well defined only
when $| s | < 1$. Therefore, one considers initial data in $H^s
\times H^{s-1}$  with
$0 < s < 1$, noticing that further smoothness would not help.   Next,
the loss of smoothness \eqref{loss}  forces us to limit $t$ to
an  interval where $0 < s  - \lambda t $, yielding only local in time
existence theorems.
We also refer to \cite{CL} for further discussions on the sharpness
of LL smoothness.

However, the local uniqueness of the Cauchy problem and the finite speed of
propagation for local solutions are  not proved in \cite{CL}.
The main goal of this paper is
to address these questions. Classical methods such as convexification,
leads one to consider general equations \eqref{maineq} with LL
coefficients in all variables.
However, the meaning of the Cauchy problem for such equations is not
completely obvious:  as mentionned above, the maximal expected
smoothness of the solutions
is $H^s$ with  $s < 1$ and their traces on the initial manifold are
not immediately
defined.  More importantly,  in the general theory of smooth
operators, the traces are defined
using
partial regularity results in the normal direction;  in our case, the
limited smoothness of the coefficients
is a source of difficulties.
It turns out that when $ s \le \mez$, one cannot   in general define
the traces of all
the first order derivatives of $u$, but only the Neumann trace
relative to the operator, which we now introduce.
Consider a smooth hypersurface $ \Sigma$ and, near a point
$ \underline y \in \Sigma $,
a  vector field conormal to $\Sigma$,  $ \nu_\Sigma  \in T^*_\Sigma$, with
$\nu_\Sigma ( \underline y) \ne 0$.
Denoting by $ \eta $ the frequency variables dual to $y$,   define
the vector field $X_\Sigma$ with symbol
\begin{equation}
\label{defXl}
X_\Sigma  (y, \eta)  = g_2 (  y ;   \eta , \nu_\sigma(y) )
\end{equation}
where $g_2 (y; \cdot, \cdot) $ is the quadratic form defined by the
principal part of $L$.

\begin{ass}
\label{assl}
$L$ is  a second order operator  of the form
$\eqref{maineq}$ on a neighborhood $\Omega$ of
$\underline y$, with coefficients $a_{j, k} \in LL (\Omega)$,
$b_j$ and $c_j$ in    $ C^\alpha(\Omega)$, for some $\alpha \in
]\frac{1}{2}, 1[$ and
$d \in L^\infty(\Omega)$.
$\Sigma $ is a smooth hypersurface through $\underline y$  and   $L$
is strictly hyperbolic    in the direction conormal to   $\Sigma$.
\end{ass}

Shrinking $\Omega$ if necessary, we assume that $ \Sigma $ is defined by
the equation $ \{ \vp = 0 \} $ with $\vp $ smooth and $d \vp \ne 0$.
We consider the one-sided Cauchy problem, say on the component
   $\Omega_+ =  \Omega \cap  \{ \vp >  0 \} $.
As usual, we say that
$u \in H_{loc}^s(\Omega \cap  \{ \vp \ge  0 \} ) $, if for any 
relatively compact open subset $\Omega_1$  of  $\Omega$, the
restriction of $u$ to
$\Omega_1  \cap  \{ \vp >   0 \} $  belongs to $H^s (\Omega \cap  \{
\vp  >  0 \} )$.

\begin{lem}
\label{traces}
i)   For all $ s  \in ] \alpha-1 ,  \alpha [ $  and $u \in
H_{loc}^s(\Omega \cap  \{ \vp \ge  0 \} ) $,
all the terms entering in the defintion of $L$ are well defined as
distributions
in  $H_{loc}^{s-2}(\Omega \cap  \{ \vp \ge  0 \})$.

ii)  If  $u \in H_{loc}^s(\Omega \cap  \{ \vp \ge  0 \}) $ and $ Lu
\in L_{loc}^2(\Omega \cap  \{ \vp \ge  0 \})$,
then the traces  $  u_{| \Sigma} $ and $( X_\Sigma  u)_{| \Sigma} $
are well defined in $H_{loc}^{ s-\mez}(\Sigma \cap \Omega)$  and
$H_{loc}^{ s-\frac{3}{2}}(\Sigma\cap \Omega)$, respectively.
\end{lem}

With this Lemma, the Cauchy problem with source term in $L^2$ and solution
in $H^s$, $s > 1 - \alpha $, makes sense.

\begin{theo} [Local existence]
\label{Th1}
Consider $s > 1- \alpha  $ and a neigborhood $\omega$ of $\underline
y$ in $\Sigma$.
Then there are $s'  \in ] 1-\alpha, \alpha [ $ and  a neighborhood
$\Omega'$ of $\underline y$  in $\RR^{1+n}$
such that  for
all Cauchy data $(u_0, u_1)$ in $H^s (\omega) \times H^{s-1}(\omega)$
near $\underline y$
and   all $f \in L^2(\Omega' \cap  \{ \vp >  0 \})$  the
Cauchy problem
\begin{equation}
\label{CP}
L u = f  , \quad  u_{| \Sigma} = u_0, \quad ( X_\Sigma  u) _{| \Sigma} = u_1,
\end{equation}
has a   solution $u \in H^{s'}(\Omega' \cap  \{ \vp >   0 \})$.
\end{theo}

   \begin{theo}[Local uniqueness]
\label{Th2}
If $s > 1 - \alpha$ and
   $u \in H^s(\Omega \cap  \{ \vp >   0 \})$ satisfies
   \begin{equation}
\label{CP0}
L u = 0 , \quad  u_{| \Sigma} =  0, \quad  (X_\Sigma  u)_{| \Sigma} = 0,
\end{equation}
then $u = 0$  on a neighborhood of $\underline y$ in $\Omega \cap  \{ 
\vp \ge  0 \}$.

   \end{theo}

   \begin{rem}
   \textup{If the coefficients of the first order term $L_1$ (see 
\eqref{defL1}) are also
$LL$, the statements above are
   true with $\alpha = 1$ since the coefficients are then  $C^\alpha$ 
for all $\alpha < 1$.
   If the $b_j$ are $C^\alpha$ and the $c_j $ are $C^{\tilde \alpha}$,
the conditions are
   $1 - \tilde \alpha < \alpha$ and the limitation  on $s$ is
   $  1 - \tilde \alpha < s$.  }
   \end{rem}

\begin{rem}
\textup{Theorem \ref{Th2} implies that if $u $ is  in $H^s$ and
satisfies $ Lu = 0$ near $\underline y$
and if $u$ vanishes on $  \{ \vp <   0 \}  $, then $u$ vanishes on a
neighborhood of $\underline y$
(see Section~5.2). Moreover, this local propagation of zero across
any space-like manifold implies
    finite speed of propagation by classical arguments which we do not
repeat here.
In particular, if $\Omega' \cap  \{ \vp \ge  0 \}$ is contained in
the domain of dependence
of $\omega$, there is existence and uniqueness for the Cauchy problem
\eqref{CP} in
    $\Omega' \cap  \{ \vp \ge  0 \}$.  }
\end{rem}

The proof of these results is given in Section 5 below.
Because all the hypotheses are invariant under smooth changes of coordinates,
we can assume that in the coordinates $y = (t, x)$, the initial surface
is $\{ t = 0\}$,
and in these coordinates, we prove the existence and uniqueness theorems.
We deduce them from similar results on   strips
$]0, T[ \times \RR^n $ and there, the  main part of the work is to
prove good energy estimates
for (weak) solutions. In this framework,
the results of Theorem~\ref{Th1} are improved, by using non isotropic spaces,
and by making a detailed account of the loss of spatial smoothness as
time evolves, as in
\cite{CGS, CL}.
The precise  results are stated in section 2 below and are proved in 
section 4 using
    the paradifferential  calculus of J.-M.  Bony, whose LL-version is
presented in section 3.

%%%%%%%%%%%%%%%%%%%%%%%%%%%%%
%     SECTION 2

\section{The global in space problem}

In this section we denote by  $ (t, x) $ the space-time 
variables.
On
$\Omega =  [0, T_0] \times \RR^n$  consider  a second order
hyperbolic differential operator
\begin{equation}
\label{eq1}
L u  = L_2 u   + L_1 u  + d  u
\end{equation}
with
\begin{eqnarray}
\label{defL2}
L_2  & = &  \D_t a_{0} \D_t   +  \sum_{j=1}^n  (  \D_t a_j \D_{x_j} +
\D_{x_j}  a_j \D_{t})
-    \sum_{j, k=1}^n \D_{x_j} a_{j, k}  \D_{x_k} ,\\
\label{defL1}
L_1 & = &   b_0 \D_t  +  \D_t c_0    +
   \sum_{j=1}^n ( b _j \D_{x_j}  +  \D_{x_j} c_j ) .
\end{eqnarray}
The coefficients satisfy on $\Omega = [0, T_0] \times \RR^n$
\begin{eqnarray}
\label{cf1} a_{j, k} = a_{k, j}  ,  \quad   a_0, \ a_j, \ a_{j, k}
\in L^\infty (\Omega) \cap  LL(\Omega) ,
\\
\label{cf2}
b_0, \ c_0 , \ b_j, \ c_j  \in L^\infty (\Omega) \cap C^\alpha(\Omega),
\\
\label{cf3}
d \in L^\infty(\Omega),
\end{eqnarray}
for some $\alpha \in ]\frac{1}{2}, 1[$. Recall that the space LL is
defined by \eqref{LL}, the semi norm
$\| a \|_{LL} $ being the best constant $C$ in $\eqref{LL}$.
In addition, for $\alpha \in ]0, 1[$, $C^\alpha$ denotes the usual
H\"older space, equipped with the norm
\begin{equation}
\label{normHol}
\| a \|_{C^\alpha} =  \| a \|_{L^\infty}   +  \sup_{ y \ne y'}
\frac{ | a(y) - a(y')|} { | y - y'|^\alpha}.
\end{equation}
When $\alpha = 1$, this defines the norm
$\| a \|_{Lip}$  in the space of Lipschitz functions.

We assume that $L$ is hyperbolic in the direction $dt$,  which means
that  there are
$\delta_0  > 0$  and $\delta_1$   such that for all
$(t, x, \xi) \in [0, T_0] \times \RR^n\times \RR^n $
\begin{equation}
\label{hyp}
a_0(t, x)  \ge  \delta_0  ,
    \quad    \sum_ { 1 \le j, k \le n} (     a_{j, k } + \frac{ a_j a_k
}{a_0}) \xi_j \xi_k    \ge \delta_1 \,  | \xi |^2 .
\end{equation}
   We denote by
   $A_{L^\infty}$,  $A_{LL} $ and $B $ constants such that for all indices
\begin{eqnarray}
    \label{defA}
\| a_0, a_j , a_{j, k} \|_{L^\infty(\Omega) }   \le A_{L^\infty}, & \quad  &
\| a_0, a_j , a_{j, k} \|_{LL(\Omega) }   \le A_{LL}, \\
\label{defB}
\| b_0, c_0, b_j , c_j  \|_{C^\alpha (\Omega) }   \le B, &\quad &  \|
d \|_{L^\infty(\Omega) } \le B.
\end{eqnarray}

\subsection{Giving sense to the Cauchy problem}
Consider the vector fields
\begin{equation}
\label{defX}
X  = a_0 \D_t  + \sum_{j=1}^n a_j \D_{x_j}   = a_0 Y .
\end{equation}
Formal   computations immediately show that the second order part of $L$ can be
written
\begin{equation}
\label{newL}
L_2  u =  Y^*   X u  - \tilde L_2 u
\end{equation}
with
\begin{equation}
\label{defY}
Y^* v  = \D_t  v + \sum_{j=1}^n \D_{x_j} (\tilde a_j v ) ,
\qquad \tilde L_2 u =
\sum_{j, k= 1}^n \D_{x_j}\big(  \tilde a_{j, k} \D_{x_k} u \big),
\end{equation}
$\tilde a_{j,k} = a_{j, k} + a_j a_k / a_0$, and  $\tilde a_j = a_j / a_0$.
Consequently, it follows that
\begin{equation}
\label{newLb}
L u = (Y^* + \tilde b_0)   (X  +   c_0  ) u -  \tilde L_2 u
+ \tilde L_1 u +  \tilde d u
\end{equation}
with
\begin{equation}
   \label{newL1}
\tilde L_1 u=   \sum_{j=1}^n \tilde b_j \D_{x_j} u +
\sum_{j = 1}^n \D_{x_j} (   \tilde c_j   u  )
\end{equation}
and
$$
\tilde b_0 =  b_0/ a_0 , \quad    \tilde b_j = b_j - \tilde b_0 a_j , \quad
    \tilde c_j  =  c_j -  \tilde a_j c_0 , \quad  \tilde d = d - c_0 \tilde c_0.
$$
The next lemma shows that these identities are rigourous under
minimal smoothness assumption on $u$.

\begin{lem}
\label{lem41} Suppose that
     $u \in H^{ \rho} (]0,T[)\times \RR^n)$ for some $\rho \in
]1-\alpha, \alpha  [$.
Then $c u$, $X u$ and $L_1 u$ belong to  $H^{ \rho-1} (]0,T[)\times
\RR^n)$. Moreover
$L_2 u$ is well defined as a distribution in  $H^{\rho-2}(]0, T[
\times \RR^n)$.

\end{lem}

\begin{proof}
$u$ and its space-time derivatives  $(\D_tu, \D_{x_j}u )$ belong to
$ H^{\rho-1}$. Following \cite{CL},
their multiplication by a bounded LL function belong to the same
space (see also Corollary~\ref{cor36}). This shows that all the
individual terms
present in the definition of  $ Xu$ belong to $ H^{\rho-1}$ and those
occurring in
$L_2 u$ and $ Y^* X u $ are well defined in $H^{\rho-2}$ in the sense
of distributions.

Next we recall that the multiplication $(b, u) \mapsto bu$ is continuous from
$ C^\alpha \times H^s$  to $H^s$ when $ | s| < \alpha $.  This
implies that the terms
$ b \D u$ and $ \D( c u)$ that occur in $L_1 u$ and $\tilde L_1 u$ belong to
$ H^{ \rho-1} $ since $\rho \in ]1-\alpha, \alpha  [$.

The last term $ d u$ is in $L^2$, thus in  $ H^{ \rho-1} $, since $c
\in L^\infty$ and $u \in L^2$.

The identity \eqref{newL} is straighforward from \eqref{defL2} since
all the algebraic computations
make sense by the preceding remarks.
\end{proof}

Next we need   partial regularity results  in time,  showing that the traces of
$ u$ and $X u $ at $ t = 0$ are well defined, as distributions, for solutions
of    $ L u = f $. This is based on the remark that this equation
is equivalent to the system
\begin{equation}
\label{syste}
\left\{\begin{aligned}
&  Y^* v + \tilde b_0   v  = \tilde L_2 u - \tilde L_1 u - \tilde d u  + f,
\\
&    Y   u + \tilde c _0 u  =   v / a_0
\end{aligned}\right.
\end{equation}
   with $ \tilde c_0 = c_0 / a_0 $. The important remark is that, for
this system, the coefficients of
   $\D_t$, both for $u$ and $v$, are equal to 1, thus smooth. Using
the notation
   $Y = \D_t + \widetilde Y$, $Y^*  = \D_t + \widetilde Y^*$, the system reads
\begin{equation}
\label{syste2}
\left\{\begin{aligned}
&  \D_t v   = -  \tilde Y^*  v   -  \tilde b_0   v  -  \tilde L_2 u -
\tilde L_1 u - \tilde d u  + f,
\\
&    \D_t u = -  \tilde Y u  +  v / a_0.
\end{aligned}\right.
\end{equation}

\begin{lem}
\label{lem42}
Suppose that  $\rho \in ]1-\alpha, \alpha  [$ and
     $u \in H^{ \rho} (]0,T[\times \RR^n)$
is such that  $L u \in L^1 ([0,T]; H^{ \rho -1} (\RR^n)) $.
Then $ u \in L^2([0,T]; H^\rho(\RR^n))$ and  $ \D_t u \in L^2([0,T];
H^{\rho-1}(\RR^n))$.
Therefore,
    $ u \in C^0 ([0,T]; H^{\rho-\mez}(\RR^n))$.

Moreover,  $ X u   \in L^2([0,T]; H^{\rho-1} (\RR^n))$   and
    $ Xu   \in C^0 ([0,T]; H^{\rho-\frac{3}{2} }(\RR^n))$.

In particular, the traces $u_{| t = 0} $ and $Xu_{| t = 0} $ are well defined
in $H^{\rho-\mez}(\RR^n)$ and $H^{\rho - \frac{3}{2}}(\RR^n)$, respectively.
\end{lem}

\begin{proof}
{\bf a)}   We use the spaces $H^{s, s'}$ of H\"ormander (\cite{H},
chapter 2), which are defined on $ \RR^ { 1+n}$ as the spaces of
temperate distributions  such that their Fourier transform $\hat u $
satisfies
$ ( 1 +  \tau ^2 + | \xi |^2)^{s/2} ( 1 +   | \xi |^2)^{s'/2} \hat u
\in L^2$. The spaces on
$ [0, T] \times \RR^n$ are defined by restriction. In particular,
$H^{0, s'} ( [0, T] \times \RR^n) =  L^2( [0, T]; H^{s'}(\RR^n))$.  Recall that
$\D_{x_j}$ maps $ H^{s, s'}$ to $ H^{s, s'-1}$ and that
\begin{equation}
\label{hypop}
u \in H^{s, s'} , \  \D_t u \in H^{s, s'-1 } \quad \Rightarrow \quad
u \in H^{s+1, s'-1}.
\end{equation}

{\bf b)}  For $ u \in H^\rho $,  the first derivatives  of $u$,  $ \tilde d u$,
as well as $\tilde  L_1u $,    $ Xu $  and $v$
   belong to
    $H^{\rho-1} =  H^{\rho-1, 0} $, as well as their multiplication by
    a LL or  $C^\alpha$ coefficient.
    Thus  $ \tilde L_2 u $ and $   \tilde  Y^* v   $ belong to
$ H^{\rho-1, -1 } $ and
\begin{equation}
\label{syste3}
\D_t  v   =   f   +   g , \quad  f = Lu \in L^1( ]0, T[; H^{\rho-1}),
\  g \in H^{\rho-1, -1 }.
\end{equation}
Let
$$
v_0 (t) = \int_0^t f(t') dt'  \in C^0 ( H^{\rho-1}).
$$
In particular, $ v \in    L^2 (]0, T[;  H^{\rho-1}) = H^{0, \rho-1} \subset
H^{\rho-1, 0}  $,
since $ \rho-1 \le  0$. Thus, $ v  - v_0 \in H^{\rho-1, 0} $ and
$ \D_t ( v - v_0) = g \in H^{\rho-1, -1 }$. By \eqref{hypop}
$ v - v_0  \in H^{\rho, -1 }  \subset  H^{0,\rho -1 }$ since $\rho \ge  0$.

Next, reasoning for fixed time and then taking $L^2$  norms we note
that the multiplication by a LL or $C^\alpha $ function maps 
$ L^2 ( ]0, T[; H^{\rho-1}) = H^{0, \rho-1} $ into itself. Thus, by the second
equation of \eqref{syste2},
$ \D_t u  = - \tilde Yu  + v/ a_0 \in  H^{0, \rho -1 } $. This finishes
the proof of the first part of the lemma.

{\bf c)}  In particular, it  implies that
$  v = Xu + b_0 u  \in H^{0, \rho-1}  $.
Thus,   $\tilde Y^* v$ and $\tilde L_2  u$  which involve  multiplication by
$C^\alpha$ or LL function, followed by  a spatial derivative, belong to
$H^{ 0, \rho-2} $.   Therefore,  the equation 
$ g \in H^{0, \rho-2}$.
Thus applying \eqref{hypop} to
$ v  - v_0  \in   H^{0,\rho -1 }$ implies that
$ v - v_0  \in   H^{1,\rho -2 }  \subset C^0 ([0,T];
H^{\rho-\frac{3}{2} }(\RR^n))$.
Since $ | \rho -\mez | < \alpha$ and
$ u \in C^0 ([0,T]; H^{\rho-\mez}(\RR^n))$,  the product $\tilde b_0
u $ belongs to
$   C^0 ([0,T]; H^{\rho-\mez}(\RR^n))$. Since $v_0$ is also in this space,
we conclude that $ Xu  \in C^0 ([0,T]; H^{\rho-\frac{3}{2} }(\RR^n))$.
\end{proof}

\begin{rem}
\textup{If $\rho > \mez$, then the multiplication by LL functions
maps   $ H^{\rho-\frac{3}{2} } $
into itself
and we can conclude that  $ \D_t u \in C^0 ([0,T];
H^{\rho-\frac{3}{2} }(\RR^n))$, as well as all the first
derivatives of $u$, so that their traces at $t = 0$ are well defined.
When $ \rho \le \mez$, the
continuity of $\D_t u$ is not clear. However, note that the trace of
$Xu $ has an intrinsic meaning,
as $X$ is the Neumann derivative on the  initial surface $\{ t = 0\}$
relative to $L$ (see the definition \eqref{defXl}). }
\end{rem}

Lemma \ref{lem42} allows us to consider  the Cauchy problem
\begin{equation}
\label{Cauchy}
L u = f, \quad  u_{| t = 0}= u_0, \quad  X u _{| t =0} = u_1,
\end{equation}
when
   $f  \in   \bigcup_{\rho > - \alpha }  L^1 ([0,T]; H^{ \rho } (\RR^n)) $ and
$u \in \bigcup_{ \rho > 1 -\alpha}  H^{ \rho} (]0,T[)\times \RR^n)$.

%%%%%%%%%%%%%%%%%%%%%%%

\subsection{The main   results}

We first state  uniqueness for the Cauchy problem:

\begin{theo}
\label{Uniq}
If  $u \in \bigcup_{ \rho > 1 -\alpha}  H^{ \rho} (]0,T[)\times
\RR^n)$ satisfies
\begin{equation}
\label{Cauchy00}
L u = 0 , \quad  u_{| t = 0}= 0 , \quad  X u _{| t =0} = 0
\end{equation}
then $u = 0$.
\end{theo}

As in \cite{CGS, CL},  we prove existence  of solutions in Sobolev spaces
having  orders decreasing in time. The proper definition is given as
follows.  The operators
\begin{equation}
\label{logD}
| D|  \quad \mathrm{and} \quad \Lambda := \mathrm{Log} ( 2 + | D | )
\end{equation}
are defined by Fourier transform, associated to the Fourier multipliers
$| \xi |$ and  $\mathrm{Log} ( 2 + | \xi |)$ respectively.

\begin{defi}
i)  $H^s(\RR^n) $ or $H^s$ denotes the usual Sobolev space on $\RR^n$.
$H^{ s + \frac{1}{2} log } $ and $ H^{ s - \frac{1}{2} log } $ denote
the spaces
$\Lambda^{- \frac{1}{2}} H^s$ and $\Lambda^{\frac{1}{2}} H^s$ respectively.

ii) Given parameters $\sigma$ and $\lambda$, we denote by
$\cC_{\sigma, \lambda}(T)$ the space of functions
$u$ such that for  all $t_0 \in [0, T ] $,  $ u \in C^0([0, t_0],
H^{\sigma- \lambda t_0})$.

iii)  $\cH_{\sigma \pm  \frac{1}{2} log , \lambda}(T)$ denotes the
spaces of  functions $u$
on $[0, T] $ with values in the space of temperate distributions in $\RR^n$
such that
\begin{equation}
\label{cH}
   (1 +  | D  |) ^{  \sigma - \lambda t}  \Lambda ^{ \pm  \frac{1}{2}}
u (t, \cdot)
   \in L^2 ([0, T] ; L^2( \RR^n)).
\end{equation}

iv)  $\cL_{\sigma  , \lambda}(T)$ denotes the space of  functions $u$
on $[0, T] $ with values in the space of temperate distributions in $\RR^n$
such that
\begin{equation}
\label{cL}
   (1 +  | D  |) ^{ \sigma - \lambda t}   u (t, \cdot)
   \in L^1 ([0, T] ; L^2( \RR^n)).
\end{equation}

\end{defi}

$\cC_{\sigma, \lambda}(T)$ is equipped with the norm
\begin{equation}
\label{normcC}
\sup_{t \in [0, T] }  \| u (t)  \|_{ H^{\sigma- \lambda t}}.
\end{equation}
The norms in $\cH_{\sigma\pm \frac{1}{2} log , \lambda}(T)$ and
$\cL_{\sigma  , \lambda}(T)$
are given by \eqref{cH} and \eqref{cL}. Equivalently,
$\cH_{\sigma \pm  \frac{1}{2} log , \lambda}(T)$ and  $\cL_{\sigma  ,
\lambda}(T)$ are the completions of
$C^\infty_0 ([0, T] \times \RR^n)$ for the norms
\begin{equation}
\label{normcH}
\| u \|_{\cH_{\sigma \pm  \frac{1}{2} log , \lambda}(T)}  =
\Big( \int_0^T \| u(t) \|^2_{ H^{\sigma - \lambda t \pm  \frac{1}{2}
log }} dt \Big)^{\frac{1}{2}}.
\end{equation}
and
\begin{equation}
\label{normcL}
\| u \|_{\cL_{\sigma   , \lambda}(T)}  =
    \int_0^T \| u(t) \| _{ H^{\sigma - \lambda t  }} dt .
\end{equation}

\begin{theo}
\label{mainthm}
Fix $\theta < \theta_1 $ in $ ]1- \alpha , \alpha[$. Then there are
$\lambda > 0$  and $K > 0$,
which depend only on the constants $A_{L^\infty}$, $A_{LL}$, $B$,
$ \delta_0 $, $\delta_1$, $\theta $ and $ \theta_1 $, given by 
\eqref{hyp}, \eqref{defA} 
and \eqref{defB} 
such that for
\begin{equation}
\label{fixTa}
T = \min \{ T_0, \frac{\theta_1 - \theta}{\lambda} \}
\end{equation}
   $ u_0 \in H^{ 1 - \theta}(\RR^n)$, $u_1 \in H^{- \theta}(\RR^n)$ and
$f = f_1 + f_2$ with $f_1 \in \cL_{- \theta, \lambda}(T)$ and $f_2\in
\cH_{- \theta - \frac{1}{2} log , \lambda}(T)$,  the Cauchy problem
$\eqref{Cauchy}$,
has a unique solution $u \in \cC_{1- \theta, \lambda} (T) \cap \cH_{1
- \theta + \frac{1}{2} log , \lambda}(T)$
with $\D_t u \in \cC_{- \theta, \lambda} (T) \cap \cH_{ - \theta +
\frac{1}{2} log , \lambda}(T)$.
Moreover, it  satisfies
\begin{equation}
\label{mainest}
\begin{aligned}
&  \sup_{ 0 \le t' \le t} \| u (t')  \|^2_{ H^{1 - \theta - \lambda t'}} +
    \sup_{ 0 \le t' \le t} \| \D_t u (t')  \|^2_{ H^{-\theta - \lambda t'}}  \\
   &  +\int_0^t \Big(  \| u (t' ) \|^2 _{  H^{1 - \theta - \lambda t' +
\frac{1}{2} log }} +
    \| \D_t u (t' ) \|^2 _{  H^{ - \theta - \lambda t' + \frac{1}{2}
log }} \Big) dt'
    \\
    \\
    & \le    K  \Big\{   \| u _0  \|^2_{ H^{1 - \theta }} +
       \|  u_1  \|^2_{ H^{-\theta }}  \\
   &  +    \big( \int_0^t   \| f_1 (t' ) \|  _{  H^{- \theta - \lambda
t'  }}  dt' \Big)^2  +
      \int_0^t   \| f_2 (t' ) \|^2 _{  H^{ - \theta - \lambda t' -
\frac{1}{2} log }}  dt' \Big\} .
    \end{aligned}
\end{equation}

\end{theo}

Note that  for $t \in [0, T]$,  $ 1 - \theta - \lambda t   \ge 1-
\theta_1 >  1 - \alpha  $, so that
    $ f \in L^1 ([0, T];  H^ {- \theta_2}  ) $  with $ \theta_1 <
\theta_2 < \alpha $. Similarly,
     $  u \in L^2 ([0, T];  H^ {1 - \theta_1}  ) $ and $  \D_t u \in
L^2 ([0, T];  H^ { - \theta_1}  ) $
implying that $  u \in  H^ {1 - \theta_1}  ([0, T] \times \RR^n) $.
Therefore, we are in a situation where we have given sense to the
Cauchy problem.

\begin{rem}
\textup{This is a local in time existence theorem since the life span
\eqref{fixTa} is limited by the choice of  $\lambda$. Thus the
dependence of $\lambda_0$ on the coefficient is of crucial importance.
In case of  Lipschitz coefficients, there is no loss of derivatives;
this would correspond to
$\lambda  = 0$. Using the notations in \eqref{defA} \eqref{defB}  and 
\eqref{hyp},
the analysis of the proof below shows that there is a function $ K_0 (
\cdot)$ such that one can  choose }
\begin{equation}\label{lda0}
\lambda   =  \frac{A_{LL}}{ \min \{ \delta_0, \delta_1\} }  \,  K_0
\big( \frac{A_{L^\infty}}{\delta_0} \big),
\end{equation}
\textup{revealing  the importance of the LL-norms of the coefficients and
the role of the hyperbolicity constant
$ \delta_1 / \delta_0$. In particular, it depends only on the second
order part of operator $L$. }
\end{rem}

\begin{rem}
\textup{A closer inspection of the proof, also shows that if the
coefficients of the pricipal part of
$L $  are $ ( a_0, a_j , a_{j, k} ) =   ( a'_0+ a''_0   , a'_j+ a''_j
, a'_{j, k} + a''_{j,k} )  $  with
$ ( a'_0, a'_j , a'_{j, k} ) $  Lipschitz continous  and $  ( a''_0,
a''_j , a''_{j, k} ) $  Log Lipschitz,
   with  LL norm bounded by  $A''_{LL}$, one can replace  $A_{LL}$ by
$A''_{LL}$ in the definition
of $\lambda$.   In particular if instead of \eqref{LL}  the
coefficients satisfy }
\begin{equation}
\label{modu}
| a( y) - a(y') |  \le C \omega ( | y - y'|)
\end{equation}
\textup{with a modulus of continuity $\omega$ such that }
\begin{equation}\label{sousLL}
\lim_{ \eps \to 0+ } \ \frac{ \omega(\eps)} { \eps | \Log \eps | }  = 0,
\end{equation}
\textup{they can be approximated  by Lipschtiz functions  with  errors
arbitrarily small  in the LL norm.
This can be done by usual mollifications, which will preseve the
$L^\infty$ bounds $A_{L^\infty}$ and keep uniform  hyperbolicity
constants $ \delta_0$ and
$\delta_1$.  As a consequence, $\lambda$ can be taken arbitrarily
small, yielding global in time
existence  with arbitrarily small loss of regularity (see Theorem~2.1
in \cite{CC} when the coefficients depend only on time).   }
\end{rem}

%%%%%%%%%%%%%%%%%%%%%%%%%%%%%%

\section{Paradifferential calculus with LL coefficients}

In this section we review several known results on paradifferential
calculus and give
the needed extensions to the case of Log-Lipschitz coefficients.

\subsection{The Paley-Littlewood analysis}

Introduce $\chi\in C^\infty_0 (\RR)$, real valued, even and such that
$0 \le \chi \le 1$ and
\begin{equation}
\label{Achi}
\chi(\xi  ) \, = 1  \quad {\rm for } \
\vert \xi  \vert \, \le 1.1 \, ,
\quad
\chi(\xi  ) \, = 0 \quad {\rm for } \
\vert \xi \vert \, \ge 1.9 \, .
\end{equation}
For $k \in \ZZ$, introduce
$\chi_k(\xi) := \chi\big(2^{-k}   \xi  \big) $,
$\widetilde \chi_k (x)$ its  inverse Fourier transform with respect
to $\xi$   and   the operators
\begin{equation}
\label{defsk}
\begin{aligned}
S_k  u \, :=\, \widetilde \chi _k * u \, =
\chi_k (D_x ) u\,,\\
\Delta_0 = S_0 ,
   \quad \mathrm{and\ for \ }  k \ge 1 \quad
\Delta _{k}   = S _{k} - S_{k-1}  .
\end{aligned}
\end{equation}
We note that $\Delta_k$ and $ S_k$ are self adjoint. Moreover, by
evenness,   $ \widetilde \chi_k $ is real, so that
$\Delta_k$ and $ S_k$ preserve reality.
For all temperate distributions $u$ one has
\begin{equation}
\label{L-P}
u =     \sum_{k \ge 0} \Delta_k u \, .
\end{equation}

The next propositions immediately follow from the definitions.

\begin{prop}
\label{propB1}
Consider $s \in \RR$.
A temperate distribution $u$ belongs to $H^s(\RR^n)$ [resp. $ H^{s
\pm \frac{1}{2}log}$]   if and only if

\quad i) for all $k \in \NN$, $\Delta_k  u \in L^2(\RR^d)$.

\quad ii)  the sequence
$\delta_k = 2^{ ks} \Vert \Delta_k  u\Vert_{L^2(\RR^d)}$
[resp.   $\delta_k =  (k+1)^{\pm\frac{1}{2}} 2^{ ks} \Vert \Delta_k
u\Vert_{L^2(\RR^d)}$]
belongs to $\ell^2(\NN)$.

Moreover,
    the norm of the sequence $\delta_k$ in $\ell^2$ is equivalent to
the norm of $u$ in the given space.

\end{prop}

\begin{prop}
\label{propB2}
Consider $s \in \RR$  and $R > 0$.  Suppose that
$\{u_k\}_{k \in \NN}$ is a sequence of functions in $L^2(\RR^d)$such that:

\quad i)    the spectrum of $ u_0$ is contained in $\{ | \xi | \le R
\}$ and for
$k \ge 1$ the spectrum of   $u_k$ is contained
in $ \left\{  {1 \over R}  \, 2^{ k} \le
| \xi |  \le  R\, 2^{k} \right\}$.

\quad ii) the sequence
$\delta_k = 2^{ ks} \Vert u_k\Vert_{L^2(\RR^d)}$
[resp.   $\delta_k =  (k+1)^{\pm\frac{1}{2}} 2^{ ks} \Vert \Delta_k
u\Vert_{L^2(\RR^d)}$]
belongs to $\ell^2(\NN)$.

Then $u = \sum u_k$ belongs to $ H^s(\RR^d)$
   [resp. $ H^{s \pm \frac{1}{2}log}$]. Moreover,
    the norm of the sequence $\delta_k$ in $\ell^2$ is equivalent to
the norm of $u$ in the given space.

When $s > 0$, it is sufficient to assume that the spectrum of
$u_k$ is contained in
$ \left\{ | \xi |   \le
R\, 2^{k} \right\}$.

\end{prop}

Next we collect several results about  the dyadic analysis of LL spaces.

\begin{prop}
\label{propB3a}
There is a constant $C$ such that for all $a \in LL(\RR^n)$ and all integers
$k  > 0$
\begin{equation}
\label{dyadLL}
\|  \Delta_k  a  \|_{L^\infty}  \le  C  k   2^{-k} \| a \|_{LL} .
\end{equation}
Moreover, for all $k \ge 0$
\begin{eqnarray}
\label{dyadLL2}
&& \| a -  S_k  a  \|_{L^\infty}  \le  C    (k+1) \| a \|_{LL}
\\
&& \label{dyadLL1}
\| S_k  a  \|_{Lip}  \le  C \Big( \| a \|_{L^\infty}  +   (k+1) \| a
\|_{LL} \Big)    .
\end{eqnarray}

If $\alpha \in ]0, 1[$ and  $a \in C^\alpha (\RR^n)$, then
\begin{equation}
\label{dyadHol}
\|  \Delta_k  a  \|_{L^\infty}  \le  C      2^{- \alpha k} \| a \|_{C^\alpha} .
\end{equation}

\end{prop}

\begin{proof} $S_k$ is a convolution operator with $\widetilde
\chi_k$ which is uniformly bounded in $L^1$. Thus
\begin{equation}
\label{dyadLinfty}
\| S_k  a  \|_{L^\infty}  \le  C   \| a \|_{L^\infty}.
\end{equation}
Moreover, since the integral of $\D_j \widetilde \chi_k$ vanishes
$$
\D_j (S_k a) (x)  = \int \D_j \widetilde \chi_k (y) \big( a (x- y) -
a(x) \big) dy .
$$
Using the LL smoothness of $a$ yields
\begin{equation}
\label{DSk}
\| \nabla S_k  a  \|_{L^\infty}  \le  C  (k+1)  \| a \|_{LL}.
\end{equation}
This implies \eqref{dyadLL1}. The proof of \eqref{dyadLL} is similar
(cf \cite{CL}).
The third  estimate is classical.
\end{proof}

%%%%%%%%%%%%%%%%%%%%%%%%%%%%%%%%%%%%

\subsection{Paraproducts}

Following J.-M. Bony (\cite{Bo}), for $N \ge 3$ one  defines the
para-product of $a$ and $u$ as
\begin{equation}
\label{parap}
T^N_a  u  = \sum_{ k = N}^\infty S_{k-N} a \   \Delta_k u
\end{equation}
   The remainder $R^N_a u$ is defined as
\begin{equation}
\label{rem}
R^N_a u  = a u - T^N_a u.
\end{equation}

   The next proposition  extends classical results (see \cite{Bo, Mey})
to the case
   of LL coefficients and Log Sobolev spaces.

\begin{prop}
\label{proppara} i) For $a \in L^\infty$ and $s \in \RR$,
$T^N_a$  continuously maps $H^s$ to $H^s$ and $H^{ s \pm \frac{1}{2} log}$
to $H^{ s \pm \frac{1}{2} log}$. Moreover, the operator   norms are
uniformly bounded
    for $s$ in a compact set.

ii)  If $a \in L^\infty \cap  LL$ and $N ' \ge N \ge 3$, $T^{N}_a -
T^{ N'}_a$ maps
$H^{ s + \frac{1}{2} log}$ into $H^{ s + 1 -  \frac{1}{2} log}$, for
all $s \in \RR$.

iii)  If $a \in L^\infty \cap LL$,  $  N \ge 3$ and $s \in ]0, 1 [$,
$R^N_a$  maps
$H^{  - s + \frac{1}{2} log}$ into $H^{ 1 - s -  \frac{1}{2} log}$, and
\begin{equation}
\label{loss1}
\Vert R^N_a u \Vert_{ H^{ 1 - s -  \frac{1}{2} log}}  \le  C
\| a \|_{LL}  \Vert  u \Vert_{ H^{- s +  \frac{1}{2} log}}
\end{equation}
with C  uniformly bounded
for $s$ in  a compact  subset of $]0, 1[$.
\end{prop}

\begin{proof}
The first statement is an immediate consequence of
\eqref{dyadLinfty} and Propositions~\ref{propB1} and  \ref{propB2}.

Next, $T^{N}_a u -  T^{ N'}_a u =  \sum_k  v_k $  with
$ v_k =    ( S_{k-N} a  - S_{k - N'} a ) \   \Delta_k u  $.   By
Proposition~\ref{propB3a}
$$
\Vert v_k \Vert_{L^2} \le   C (k+1) 2^{-k} \Vert \Delta_k u \Vert_{L^2} .
$$
With Proposition~\ref{propB2}, this implies $ii)$.

To prove $iii)$ we can assume that $N =3$. Then
\begin{equation}
\label{reste}
R_a u = \sum_{k \ge 3} \Delta_k a\ S_{k-3}  u   +
\sum_k  \sum_{| k - j |  \le 2 }  \Delta_j a \Delta_k u .
\end{equation}
If $u \in  H^{ - s +   \frac{1}{2} log}$, then
\begin{equation*}
\|  \Delta_j u \|_{L^2}  \le
\frac{2 ^{j s} }{\sqrt{j+1}}   \eps_j
\end{equation*} 
with $\{ \eps_j\} \in \ell^2$. We note that the sequence
\begin{equation}
\label{resomm}
\widetilde \eps_k  = \sum_{ j \le k}  \frac{\sqrt{k+1}  }{\sqrt{j+1}}
2 ^{(j- k) s } \eps_j
\end{equation}
is also in $ \ell^2$ with
\begin{equation*}
\| \widetilde \eps_k \|_{\ell^2}  \le  C
\|  \eps_j  \|_{\ell^2}
\end{equation*}
with $C$ uniformly bounded when $s $ in a compact subset of $]0 , +
\infty ]$. Thus
\begin{equation*}
\|  S_{k-3} u \|_{L^2}  \le
\frac{2 ^{k s} }{\sqrt{k+1}}     \eps'_k
\end{equation*}
with $\{\eps'_k\} \in \ell^2$.
Therefore,
\begin{equation*}
   \|   \Delta_k a\ S_{k-3}  u   \|_{L^2}  \le
C  \sqrt{k+1}\    2^{ (s-1) k }     \eps'_k .
\end{equation*}
Proposition~\ref{propB2} implies that the first sum in \eqref{reste} belongs to
$  H^{ 1 - s -   \frac{1}{2} log}$.

Similarly,
\begin{equation*}
\Big \|    \sum_{| k - j |  \le 2 }  \Delta_j a \Delta_k u    \Big
\|_{L^2}  \le
C  \sqrt{k+1}\    2^{ (s-1) k }     \eps''_k .
\end{equation*}
   with $\{\eps''_k\} \in \ell^2$. Now the spectrum of  $\Delta_j a
\Delta_k u  $ is
   contained in the ball $\{ | \xi | \le  2^{k +3}\}$; because $1- s > 0$,
   Proposition~\ref{propB2} implies that the second sum in
\eqref{reste} also belongs to
$  H^{ 1 - s -   \frac{1}{2} log}$, and the norm is uniformly bounded when
$s$ remains in a compact subset of $[0, 1[$.
\end{proof}

\begin{rem}
\textup{By $ii)$ we see that the choice of $N \ge 3$ is essentially
irrelevant in our analysis, as in
\cite{Bo}. To simplify notation, we make a definite choice of $N$,
for instance
$N = 3$, and use the notation $T_a$ and $R_a$ for $T^N_a$ and $R^N_a$. }

\end{rem}

\begin{cor}
\label{cor36}
The multiplication $ (a, u) \mapsto au$ is continuous from $ (
L^\infty\cap LL) \times
H^{s + \delta log} $  to $ H^{s + \delta log} $ for $ s \in ]-1, 1[$
and $ \delta \in \{-\frac{1}{2}, 0,  \frac{1}{2} \}$.
\end{cor}

\begin{proof}
(see \cite{CL}) Property $iii)$ says that $ R_a$ is smoothing by
almost one derivative in negative
spaces, and therefore, for all
$ \sigma \in ]-1, 1[$  it  maps $ H^ \sigma$ to $H^{\sigma'} $
for all  $ \sigma' > \max{ \{\sigma, 0}\} $  such that $ \sigma' <
\min{ \{\sigma + 1 , 1\} }$.
Combining this observation with $i)$, the corollary follows.
\end{proof}
In particular, we note  the following estimate
\begin{equation}
\label{estprod}
\| a u \|_{ H^{s + \frac{1}{2}  log}}  \le  C \big(  \| a
\|_{L^\infty} \|  u \|_{ H^{s + \frac{1}{2}  log}}   +
   \| a \|_{LL} \|  u \|_{ H^{s }}  \big).
\end{equation}

\begin{prop}
\label{propB3b}   Consider $q = \sqrt{ (1 +  |\xi | ^2) } $ and
$\psi (\xi)$  a symbol of degree $m$   on $\RR^n$. Denote by
    $Q  = \sqrt {(1 - \Delta)}  $  and $\Psi$ the associated operators.
    If $a \in L^\infty \cap LL$, then
    the commutator $[ Q^{- s} \Psi    , T_ a ] $ maps
   $ H^{  - s +  \frac{1}{2} log} $ into $ H^{ 1-m-  \frac{1}{2} log} $ and
   \begin{equation}
   \label{loss2}
\Vert [ Q^{- s} \Psi   , T_ a ]  u \Vert_{ H^{ 1-m   -  \frac{1}{2}
log}}  \le  C
\| a \|_{LL}  \ \Vert  u \Vert_{ H^{- s +  \frac{1}{2} log}}
\end{equation}
with C  uniformly bounded
for $s \in [ 0, 1] $ and $\psi$ in a bounded set.
\end{prop}

\begin{proof}
   We use Theorem 35 of \cite{CM}, which states that if $ H $ is a
Fourier multiplier with
symbol  $h$ of degree $0$ and if $a$ is Lipschitzean, then
$$
\| [ H, a ]  \D_{x_j} u \|_{L^2} \le  C \| \na_x a \|_{L^\infty}  \
\| u \|_{L^2}.
$$
For $k > 0$, writing $ \Delta_k u$ as sum of derivatives, this implies that
\begin{equation}\label{comCM}
\| [ H, a ]  \Delta_k  u \|_{L^2} \le  C  2^{- k}  \| \na_x a
\|_{L^\infty}  \ \| \Delta_k u \|_{L^2}.
\end{equation}
with $C$ independent of $k$ and $H$, provided that the symbol $h$
remains in a bounded
set of symbols of degree $0$.

We now proceed to the proof of the proposition.
Since $\Psi$ and $Q$ commute with $\Delta_k$, one has
\begin{equation}
\label{commut}
[ Q^{- s} \Psi  , T_ a ]  u =  \sum_{ k \ge 3}  [ Q^{- s} \Psi ,
S_{k-3} a  ]  \Delta_k  u .
\end{equation}
Moreover,  since the spectrum of $S_{k -3} a \Delta_k u$  is contained
in the annulus $ 2^{ k-1} \le | \xi | \le 2^{k + 2}$, it follows that
\begin{equation}
\label{red}
   [ Q^{- s} \Psi , S_{k-3} a  ]  \Delta_k  =   2^{  k (m-s) }  [  H_k
, S_{k-3} a  ]  \Delta_k
\end{equation}
where   the symbol of $H_k$ is
$$
h_k ( \xi)  = 2^{k (s-m) } q ^{-s}(\xi) \psi( \xi) \vp( 2^{ - k} \xi )
$$
and $\vp$ supported in a  suitable fixed annulus. Note that the family
   $\{ h_k\} $  is bounded in the space of  symbols of degree $0$,
uniformly in $k$, $s \in [0, 1]$ and $\psi$ in a bounded set of
symbols of degree $m$.
   By \eqref{comCM}, it follows that
\begin{equation*}
\Vert  [  H_k   , S_{k-3} a  ]   \Delta_k u   \Vert _{L^2} \le
C  2^{ k (m-s-1)}  \Vert \na S_{k-3} a \Vert_{L^\infty}  \ \Vert
\Delta_k u \Vert_{L^2}.
\end{equation*}
Together with \eqref{DSk} and Proposition~\ref{propB1}, this implies that for
$u \in H^{ -s + \frac{1}{2} log}$,
\begin{equation*}
\Vert
   [ Q^{- s} \Psi  , S_{k-3} a  ]  \Delta_k  u \Vert \le C  (k+1)
\Vert a \Vert_{LL} \
   \Vert \Delta_k u \Vert_{L^2}.
\end{equation*}
    Using  Proposition~\ref{propB2}, the estimate
\eqref{loss2} follows.
\end{proof}

\begin{prop}
\label{prop37}
If $a \in L^\infty \cap LL$ is real valued, then
      $  \big( T_ a  - (T_a)^*\big) \D_{x_j}  $   and  $ \D_{x_j}
\big( T_ a  - (T_a)^*\big)   $ map
   $ H^{  0 +  \frac{1}{2} log} $ into $ H^{ 0 -  \frac{1}{2} log} $ and satisfy
   \begin{equation}
   \label{loss3}
\begin{aligned}
\Vert  \big( T_ a  - (T_a)^*\big) \D_{x_j}    u \Vert_{ H^{ 0 -
\frac{1}{2} log}}  \le  C
\| a \|_{LL}  \ \Vert  u \Vert_{ H^{0 +  \frac{1}{2} log}} ,
\\
\Vert \D_{x_j}   \big( T_ a  - (T_a)^*\big)   u \Vert_{ H^{ 0 -
\frac{1}{2} log}}  \le  C
\| a \|_{LL}  \ \Vert  u \Vert_{ H^{0 +  \frac{1}{2} log}} .
\end{aligned}
\end{equation}
   \end{prop}

\begin{proof} The  $S_k a$ are  real valued, since $a$ is real, and
the $\Delta_k$ are
self adjoint, thus
\begin{equation*}
(T_a)^*  u  =   \sum_{k=3}^\infty
    \Delta_k \big( (   S_{k-3} a ) \,u   \big) .
\end{equation*}
Therefore, one has
\begin{equation*}
   \big( T_ a  - (T_a)^*\big)   =  \sum
   [ S_{k-3} a, \Delta_k]   =   \sum
   [ S_{k-3} a, \Delta_k]  \Psi_k
\end{equation*}
where
$\Psi_k$ is a Fourier multiplier with symbol $\psi_k =\psi (2^{-k} \xi)$ and
$\psi$ is supported in a suitable annulus. Using again \cite{CM} (see
\eqref{comCM}) yields
\begin{equation*}
\label{318}
   \Vert  [ S_{k-3} a, \Delta_k] \D_{x_j} \Psi_k  u \Vert_{L^2}
   \le C (k+1) \Vert a \Vert_{LL} \Vert \Psi_k u \Vert_{L^2},
\end{equation*}
and a similar estimate when the derivative is on the left of the commutator.
Since the spectrum of $ [ S_{k-3} a, \Delta_k]   \Psi_k  u$ is contained in
a annulus of size $\approx 2^k$, this implies \eqref{loss3}.
\end{proof}

\begin{prop}
\label{prop38}
If $a$ and $b$ belong to  $ L^\infty \cap LL$, then
      $  \big( T_ a T_b  - T_{ a b } \big) \D_{x_j}  $ maps
   $ H^{  0 +  \frac{1}{2} log} $ into $ H^{ 0 -  \frac{1}{2} log} $ and
   \begin{equation}
   \label{loss4}
   \begin{aligned}
\Vert   \big( T_ a T_b  - T_{ a b } \big) \D_{x_j}  &  u \Vert_{ H^{
0 -  \frac{1}{2} log}}
\\
   &  \le  C
\Big( \| a \|_{LL}  \| b \|_{L^\infty} +   \| b \|_{LL}  \| a
\|_{L^\infty}  \Big)
\ \Vert  u \Vert_{ H^{0 +  \frac{1}{2} log}} .
\end{aligned}
\end{equation}

\end{prop}

\begin{proof}
By Proposition \ref{proppara}, it is sufficient to prove the estimate
for any paraproduct
$T^N$. One has
\begin{equation*}
T^N_a T^N_b \D_{x_j} u  =  \sum_{ k \ge N} \sum_{ l \ge N} S_{k-N} a\
\Delta_k \Big( S_{l- N} b \  \Delta_l \D_{x_j} u\Big).
\end{equation*}
In this sum,   terms with $ | l - k | \le 2$  vanish, because of the
spectral localization of
$ S_{l- N} b \,  \Delta_l \D_{x_j} $.  The commutators
$ [  \Delta_k,   S_{l- N} b ]  $ contribute to terms  which are
estimated as in \eqref{commut}:
\begin{equation*}
\Vert  [  \Delta_k,   S_{l- N} b ]   \Delta_l \D_{x_j} u \Vert_{L^2}
\le  C (k+1)  \|  b \|_{LL}  \ \Vert \Delta_l  u  \Vert_{L^2}.
\end{equation*}
If $N$ is large enough, the spectrum
of the corresponding term is contained in a annulus of size
$\approx 2^k$ and hence the commutators contribute to an error term
in \eqref{loss4}.  Therefore, it is sufficient to
estimate
\begin{equation}
\label{320}
     \sum_{ k \ge N} \sum_{ l \ge N} \Big( S_{k-N} a
   S_{l- N} b   -    S_{k- N} (ab) \Big) \Delta_k \,  \Delta_l \D_{x_j} u.
\end{equation}
Again, only terms with $ | l - k | \le 2$ contribute to the sum.
Using \eqref{dyadLL2}, one has
\begin{equation*}
\begin{aligned}
   \| a  -   S_{ k - N}  a \|_{L^\infty} &\le  C (k+1) 2^{-k} \| a \|_{LL},
   \\
    \| b  -   S_{ l  - N}  b \|_{L^\infty} &\le  C (k+1) 2^{-k} \| b \|_{LL},
    \\
     \| ab  -   S_{ k - N} ( a b)  \|_{L^\infty}& \le  C (k+1) 2^{-k}
\| a b \|_{LL}.
\end{aligned}
\end{equation*}
Thus
\begin{equation*}
\begin{aligned}
\|  S_{k-N} a
   S_{l- N} b   -   & S_{k- N} (ab)  \|_{L^\infty}  \\
   & \le
   C  (k+1) 2^{-k}   \big(  \| a  \|_{LL} \| b \|_{L^\infty}  +  \|
a\|_{L^\infty} \| b \|_{LL} \big) .
   \end{aligned}
\end{equation*}
Since the terms in the sum \eqref{320} have their spectrum in annuli
of size $\approx 2^k$, this implies that this sum belongs to $H^{0-
\frac{1}{2} log}$ when
   $u \in H^{0+ \frac{1}{2} log}$, with an estimate similar to \eqref{loss4}.
\end{proof}

%%%%%%%%%%%%%%%%%%%%%%%

\subsection{Positivity estimates}

The paradifferential calculus sketched above is
well adapted to the analysis of high frequencies but does not take into account
the low frequencies.  For instance, the positivity of the function
$a$ does not imply
the positivity of the operator $T_a$ in $L^2$, only the positivity up
to a smoothing
operator.  However, in the derivation of energy estimates, such
positivity results
are absolutely necessary.
To avoid a separate treatment of low frequencies, we introduce
{\sl modified paraproducts}  for which positivity results  hold
(we could also introcude weighted paraproducts as in   \cite{Me, Me2, MZ}).

Consider a nonnegative   integer  $\nu$ and define the modified paraproducts
\begin{equation}
\label{modparap}
P^\nu_a  u  = \sum_{ k = 0}^\infty S_{\max\{\nu, k-3\}} a \
\Delta_k u =  S_\nu a S_{\nu+2} u  +
\sum_{ k = \nu}^\infty S_{ k} a \   \Delta_{k+3} u .
\end{equation}
   Then
\begin{equation}
\label{difpara}
P^\nu_a u  - T_ a u  =   \sum_{ k = 0}^{\nu+2}
\sum_{j = \max\{0, k-2 \}}^{\nu} \Delta_j a \   \Delta_k u
\end{equation}
and
\begin{equation}
\label{difpara2}
a u -  P^\nu_a u     =   \sum_{ j = \nu +1}^\infty
    \Delta_j a \   S_{j+2} u .
\end{equation}
   The difference \eqref{difpara} concerns only low frequencies, and
therefore the results of Propositions~$\ref{propB3b}, \ref{prop37}$
and ~$\ref{prop38}$ are valid if one substitutes
$P^\nu_a$ in place of $T_a$, at the cost of additional error terms.
In particular,  \eqref{difpara} and \eqref{difpara2} immediately
imply the following estimates:

\begin{lem}
\label{lem312a}
i)      There is a constant $C$ such that for all $\nu$,
   $a\in L^ \infty$ and all $u \in L^2$,
\begin{equation}
\label{331a}
\|  ( P^\nu_a    - T_ a) \D_{x_j}  u   \|_{L^2}
+
\|   \D_{x_j} ( P^\nu_a    - T_ a)  u   \|_{L^2}
\le C  2^\nu \| a \|_{L^\infty} \, \| u \|_{L^2}.
\end{equation}

ii)   There is a constant $ C_0 $ such that for all $\nu$
for all $a\in LL$ and all $u \in L^2$,
\begin{equation}
\label{332a}
\|  a u  -  P^\nu_a u      \|_{L^2 } \le C_0  \nu 2^{-\nu}   \| a
\|_{LL} \, \| u \|_{L^2}.
\end{equation}

\end{lem}

We will also use the following extension of Proposition~\ref{prop37}:

\begin{prop}
\label{prop311b}
If $a \in L^\infty \cap LL$ is real valued, then
      $  \big( P^\nu_ a  - (P^\nu_a)^*\big) \D_{x_j}  $  and
$  \D_{x_j}  \big( P^\nu_ a  - (P^\nu_a)^*\big)  $ map
   $ H^{  0 +  \frac{1}{2} log} $ into $ H^{ 0 -  \frac{1}{2} log} $ and
   \begin{equation}
   \label{loss3b}
\begin{aligned}
&\Vert  \big( P^\nu_ a  - (P^\nu_a)^*\big) \D_{x_j}    u \Vert_{ H^{
0 -  \frac{1}{2} log}}  \le  C
\| a \|_{LL}  \big(  \Vert  u \Vert_{ H^{0 +  \frac{1}{2} log}}  +
\nu  \| u\|_{L^2} \big),
\\
& \Vert \D_{x_j}  \big( P^\nu_ a  - (P^\nu_a)^*\big)    u \Vert_{ H^{
0 -  \frac{1}{2} log}}  \le  C
\| a \|_{LL}  \big(  \Vert  u \Vert_{ H^{0 +  \frac{1}{2} log}}  +
\nu  \| u\|_{L^2} \big).
\end{aligned}
\end{equation}
\end{prop}

\begin{proof} One has
\begin{equation*}
   \big( P^\nu_ a  - (P^\nu_a)^*\big) \D_{x_j}  u   =  [ S_{\nu} a,
S_{\nu+2} ] \D_{x_j}  u +
   \sum_{k \ge \nu}
   [ S_{k} a, \Delta_{k+3} ] \D_{x_j}  u .
\end{equation*}
The sum over $k$ is  treated exactly as in the proof of
Proposition~\ref{prop37} and contibutes to the same error term.
Using again  Theorem~35 of \cite{CM}, the  $L^2$   norm of the
first  term is estimated by
$$
C  \| \nabla_x S_{\nu} a \|_{L^\infty}  \|   u \|_{L^2}
\le  C (\nu +1) \| a\|_{LL}  \| u \|_{L^2}
$$
and contibutes to the second error term in \eqref{loss3b}.  When the
derivative is on the left,
the proof is similar.
\end{proof}

Moreover, a comparison of $T_a u $ with $au$ immediately implies the
following positivity estimate.

\begin{cor}
\label{cor314a}
There is a constant $ c_0 $, such that  for any  positive LL-function
$a$ such that
$\delta =  \min a(x)   > 0 $,   all  $\nu $  such that
$ \nu 2^{- \nu}  \le c_0   \delta / \| a\|_{LL} $,
    and   $ u \in L^2(\RR^n)$,
\begin{equation}
\label{333a}
{\rm Re} \big( P^\nu_a u ,  u  \big)_{L^2}  \ge \frac{\delta}{2}  \|
u\|^2_{L^2}.
\end{equation}

\end{cor}
Here, $( \cdot  ,  \cdot)_{L^2} $ denotes the scalar product in  $L^2$.
This estimate extends to vector valued functions $u$ and matrices
$a$, provided that
$a$ is symmetric and positive.

%%%%%%%%%%%%%%%%%%%%%%

\subsection{The time dependent case}

In the sequel we will consider functions of $(t, x) \in [0, T] \times
\RR^n$, considered as functions
of $t$ with values in various  spaces of functions of $x$. In
particular we denote by
$T_a$ the operator acting for each fixed $t$ as $T_{a(t)}$ :
\begin{equation}
\label{parapb}
( T_a  u)(t)   = \sum_{ k = 3}^\infty S_{k-3}(D_x)  a(t)  \
\Delta_k (D_x)  u(t) .
\end{equation}
The Propositions~\ref{proppara}, \ref{propB3b}, \ref{prop37} and \ref{prop38}
apply for each fixed $t$.  There are similar definitions for the
modified paraproducts
$ P^\nu_a$; further, Lemma~\ref{lem312a} and Corollary~\ref{cor314a} apply for
fixed $t$.
\medbreak

When $a$ is a Lipschitz function of $t$, the definition
\eqref{parapb} immediately implies that
\begin{equation}
\label{comL}
[ \D_t ,  T_a ]  = T_{\D_t a}, \qquad [ \D_t ,  P^\nu_a ]  = P^\nu_{\D_t a} .
\end{equation}
When $a$ is only Log Lipschitz this formula does not make sense,  since
$\D_t a$  is not defined as a function.  The idea, already used in
\cite{CGS, CL}, is that
it is sufficient to  commute $\D_t$ with  time regularization of $a$.
In our context, this simply means that  in \eqref{parapb}, we will
replace the term $ S_{k-3} a $, which is a
spatial regularization of $a$,  by a
space-time  regularization, namely
$ S_{k-3} a_k $  where $a_k$ is a suitable time mollification of $a$.
We now give the details for $P^\nu$, as we will need them in the next section.

Introduce the mollifiers
\begin{equation}
\label{molli}
\jmath_{k} (t)  =  2^k  \jmath( 2^k t)
\end{equation}
where $\jmath \in C^\infty_0(\RR)$ is non negative, with integral
over $\RR$ equal to $1$.

\begin{defi}
Given $a \in L^\infty \cap LL ([0, T_0] \times \RR^n)$, define
\begin{equation}
\label{rega}
a_k(t,x)  =  \jmath_k *_t \tilde a = \int \jmath_k (t - s) \tilde a(s, x) ds
\end{equation}
where $\tilde a $ is the LL extension of $a$ given by
\begin{equation}
\label{exta}
\tilde a(t, x) = a(0, x), \quad t \le 0, \qquad \tilde a(t, x) =
a(T_0, x), \quad t \ge T_0.
\end{equation}
Next, for fixed $t$, the operator $\widetilde T_{a(t)} $ is defined by
\begin{equation}
\label{parapb2}
    \widetilde P^\nu_{a(t)}   u    = S_\nu a_\nu  S_{\nu+2} u  +
\sum_{ k = \nu }^\infty S_{ k} a_{k}  \   \Delta_{k+3} u .
\end{equation}
We denote by $\widetilde P^\nu_{a} $  the operator acting on
functions of $(t,x)$ by
$(\widetilde P^nu_a u)(t) =  \widetilde P^\nu_{a(t)}   u(t)$.
\end{defi}

\begin{prop}

\label{prop314c}
Let   $a \in L^\infty \cap LL ([0, T_0] \times \RR^n)$. Then
    for each $t \in [0,  T_0]$, the operators $R_1(t) = (P^\nu_{a(t)} -
\widetilde P^\nu_{a(t)} ) \D_{x_j} $,
$R_2(t) =   \D_{x_j}  (P^\nu_{a(t)} - \widetilde P^\nu_{a(t)} )$,
$R_3(t)  = \big( (\widetilde P^\nu_{a(t)})^*  - \widetilde
P^\nu_{a(t)} \big) \D_{x_j} $,
$R_4(t)  = \D_{x_j}  \big( (\widetilde P^\nu_{a(t)})^*  - \widetilde
P^\nu_{a(t)} ) $,  and
$R_5(t) =  [D_t , \widetilde P^\nu_{a} ] (t) $
map
   $    H^{  0 +  \frac{1}{2} log}   $ into $
    H^{ 0 -  \frac{1}{2} log} $ and there is a constant $C$ such that
for all $t \in [0,  T_0]$ and for $k=1,\dots,5$,
   \begin{equation}
   \label{loss5}
   \begin{aligned}
\Vert    R_k    u \Vert_{ H^{ 0 -  \frac{1}{2} log}}
     \le  C
    \| a \|_{LL}  \big(
\Vert  u \Vert_{ H^{0 +  \frac{1}{2} log}}   + \nu \|  u \|_{L^2} \big).
\end{aligned}
\end{equation}

\end{prop}

\begin{proof}  {\bf a)} First, we recall several estimates from \cite{CL}.
For  $a \in LL ([0, T_0] \times \RR^n)$ the difference $a- a_k$ satisfies
\begin{eqnarray}
\label{appa}
&& |   a (t, x) - a_k(t, x)  |  \le  C (k+1) 2^{ - k} \| a \|_{LL},
\\
\label{appdta}
&& |   \D_t a_k(t, x)  |  \le  C (k+1)   \| a \|_{LL}.
\end{eqnarray}
with $C$ independent of $t$ and $x$. In particular, we note that
\begin{equation}
\label{331b}
   \|  S_{ k } (a(t) - a_k(t) ) \|_{L^\infty}  \le    C (k+1) 2^{ - k}
\| a \|_{LL}.
\end{equation}

{\bf b)}
In accordance with \eqref{parapb2}, for $ l = 1,2, 5$, we split  
$R_l$ in two terms
\begin{equation}
\label{337b}
R_l (t) u =   B_l u  +  H_l u , \quad   H_l u   =   \sum_{k \ge  \nu} w_k
\end{equation}
with $ B_l u $ spectrally supported in the  ball of radius $ 2^{ \nu
+4  }$  and
with   $w_k$ spectrally supported in
an annulus $ | \xi | \approx 2^k$. For  $ R_1$,
\begin{equation*}
\label{parapbb1}
   B_1 u  =  S_{\nu} (a(t) - a_\nu (t))   \   S_{\nu+2}  \D_{x_j} u,
\quad w_k  =  S_{k} (a(t) - a_k (t))   \   \Delta_{k+3}  \D_{x_j} u  .
\end{equation*}
With \eqref{331b}, this implies that
$$
\| B_1 u  \|_{L^2} \le  C (\nu+1)  \| a\|_{LL}  \|   u \|_{L^2}
$$
and
$$
\| w_k \|_{L^2} \le  C (k+1) \| a\|_{LL}  \| \Delta_{k+3} u \|_{L^2},
$$
implying that
$$
\Vert    H_1     u \Vert_{ H^{ 0 -  \frac{1}{2} log}}
     \le  C
    \| a \|_{LL}
\Vert  u \Vert_{ H^{0 +  \frac{1}{2} log}}.
$$

\medbreak

For  $ R_2$,  the analysis is similar. One has
\begin{equation*}
\label{parapbb2}
   B_2 u =   \D_{x_j} \big( S_{\nu} (a(t) - a_\nu (t))   \   S_{\nu+2}  u \big),
\quad w_k  =   \D_{x_j} \big( S_{k} (a(t) - a_k (t))   \
\Delta_{k+3}   u  \big).
\end{equation*}
   Thanks to the spectral localization, the estimates for $B_2u $ and $w_k$ are
the same as in the case of $R_1$,
   implying that
\begin{eqnarray}
&& \| B_2 u  \|_{L^2} \le  C (\nu+1)  \| a\|_{LL}  \|   u \|_{L^2}
\\
&& \Vert    H_2     u \Vert_{ H^{ 0 -  \frac{1}{2} log}}
     \le  C
    \| a \|_{LL}
\Vert  u \Vert_{ H^{0 +  \frac{1}{2} log}}.
\end{eqnarray}

\medbreak
{\bf c)} For $k = 5$ we write \eqref{337b}  with 
\begin{equation*}
B_5 u  =  S_{\nu }   (\D_t a_\nu (t))   \   \Delta_{\nu+2}     u,
\quad
w_k = S_{k}   (\D_t a_k (t))   \   \Delta_{k+3}     u  .
\end{equation*}
Thus the estimates  \eqref{appdta}  imply
\begin{eqnarray*}
&& \| B_5 u  \|_{L^2} \le  C (\nu+1)  \| a\|_{LL}  \|   u \|_{L^2}
\\
&& \Vert    H_5     u \Vert_{ H^{ 0 -  \frac{1}{2} log}}
     \le  C
    \| a \|_{LL}
\Vert  u \Vert_{ H^{0 +  \frac{1}{2} log}}.
\end{eqnarray*}
\medbreak

{\bf c)} One has
$$
R_3  (t) =  R_1(t)  + R_2^*(t)  +  \big( ( P^\nu_{a(t)})^*   -
P^\nu_{a(t)} \big) \D_{x_j} .
$$
The third term is estimated in  Proposition~\ref{prop311b}.
The operators $R_1$  and $R_2^* =  B_2^* + H_2^* $ are estimated in part b),
implying that $R_3$ satisfies   \eqref{loss5} for $k=3$. 
The proof for
$ R_4 = R_3^* = R^*_1 + R_2  +  \D_{x_j}  \big( ( P^\nu_{a(t)})^*   -
P^\nu_{a(t)} \big) $
is similar.

This finishes the proof of the Proposition.
\end{proof}

\bigbreak

\begin{lem}
\label{lem315b}
     There is a constant $ C_0 $ such that
for any $a\in LL([0,T_0] \times \RR^n) $, $u \in L^2(\RR^n)$,
$\nu\ge 0$  and all
$ t \in [0, T_0]$, one has
\begin{equation}
\label{342b}
\|  a (t)  u  -  \widetilde P^\nu_{a(t)}  u      \|_{L^2 } \le C_0
\nu 2^{-\nu}   \| a \|_{LL} \, \| u \|_{L^2}.
\end{equation}

\end{lem}

\begin{proof}
We have
$$
   a u  -  \widetilde P^\nu_a u   =
    (a -  S_\nu a_\nu  ) S_{\nu+2} u  +
\sum_{ k = \nu }^\infty (a - S_{ k} a_{k} ) \   \Delta_{k+3} u .
$$
Combining \eqref{dyadLL2} and \eqref{331b}, we see that
$$
\|   a (t)  - S_{ k} a_{k} (t) \|_{L^\infty}  \le C k 2^{- k}  \| a \|_{LL}.
$$
This implies \eqref{342b}.
\end{proof}

The  lemma immediately implies the following positivity estimate.

\begin{cor}
\label{cor316b}
There is a constant $ c_0 $, such that  for any  positive LL-function
$a$ such that
$\delta =  \min a(t, x)   > 0 $,   all  $\nu $  such that
$ \nu 2^{- \nu}  \le c_0   \delta / \| a\|_{LL} $,
    and   $ u \in L^2(\RR^n)$,
\begin{equation}
\label{333aa}
{\rm Re} \big( P^\nu_{a(t)}  u ,  u  \big)_{L^2(\RR^n) }  \ge
\frac{\delta}{2}  \| u\|^2_{L^2(\RR^n)}.
\end{equation}
The same result holds for vector valued functions $u$ and definite
positive square matrices
$a$.

\end{cor}

Finally, we quote the following commutation result which will be
needed in the next section.

\begin{prop}
\label{prop317}
Suppose that  $a\in LL([0,T_0] \times \RR^n) $. Then $ \Lambda^\mez
[ \widetilde P^\nu_{a(t)},  \Lambda^\mez] $   and    $   [ \widetilde
P^\nu_{a(t)} ,  \Lambda^\mez]   \Lambda^\mez$
are bounded in  $L^2$ and satisfy
\begin{equation*}
\begin{aligned}
\|   \Lambda^\mez  [ \widetilde P^\nu_{a(t)},  \Lambda^\mez]  u \|_{L^2}
& + \|   [ \widetilde P^\nu_{a(t)},  \Lambda^\mez]   \Lambda^\mez  u \|_{L^2}
\\
& \le  C  \big( \nu^2 2^{-\nu} \| a \|_{LL}  + \nu \| a \|_{L^\infty}
\big)  \| u\|_{L^2}.
\end{aligned}
\end{equation*}

\end{prop}

\begin{proof}
Thanks to the spectral localisation, the low frequency part
$ S_\nu a_\nu S_{\nu + 2} $  in $\widetilde P^\nu_a$
contributes to terms whose $L^2$ norm is bounded by
$$
C  \nu  \| u \|_{L^2}.
$$
The commutator  with the high frequency part reads
\begin{equation*}
      \sum_{ k \ge \nu}  [ \Lambda^\mez  , S_{k} a_k  ]  \Delta_{k+ 3}  u .
\end{equation*}
We argue as in the proof of Proposition~\ref{propB3b} and write
\begin{equation}
\label{redb}
[ \Lambda^\mez  , S_{k} a_k  ]  \Delta_{k+ 3}   =   ( k+1)^\mez\,  [
H_k   , S_{k} a_k  ]  \Delta_{k+3}
\end{equation}
where   the symbol of $H_k$ is
$
h_k ( \xi)  =  (k+1)^{-\mez}  ( \Log ( 2 + | \xi |) )^\mez\vp( 2^{ - k} \xi )
$
and $\vp$ is supported in a  suitable fixed annulus. Note that the family
   $\{ h_k\} $  is bounded in the space of  symbols of degree $0$.
    By \eqref{comCM}, one has
\begin{equation*}
\Vert  [  H_k   , S_{k} a (t)  ]   \Lambda^\mez \Delta_{k+3}  u
\Vert_{L^2}   \le
C    ( k+1)    2^{-k}   \|  \na_x S_k a_k(t) \|_{L^\infty}  \Vert
\Delta_k u \Vert_{L^2}.
\end{equation*}
Since  $ \na_x S_k a_k   =  (  \na_x S_k a) * \jmath_k$, its $L^\infty$ norm
is bounded by $ C k \| a \|_{LL}$.
Adding up, and using the spectral localization, these terms contribute
a function  whose $L^2$ norm is bounded by
$ C \nu^2 2^{-\nu} \| a \|_{LL}     \| u\|_{L^2} $.

When $\Lambda^\mez$ is on the left of the commutator, the analysis is similar.
\end{proof}

%%%%%%%%%%%%%%%%%%%%%%%%%%%%%
%%%%%%%%%%%%%%%%%%%%%%%%%%%%%

%   SECTION   4

%%%%%%%%%%%%%%%%%%%%

\section{Proof of  the main results}

%%%%%%%%%%%%%%%%%%%%%%%%%%%%%

\subsection{The main estimate}

We consider the operator \eqref{eq1} with coefficients which satisfy
\eqref{cf1}, \eqref{cf2} and
\eqref{cf3}.
We fix $\theta  < \theta_1$  in $ ]1- \alpha, \alpha[$, and with
$\lambda$  to be chosen later,
we introduce the notation
\begin{equation}
\label{st}
s(t) = \theta  + t \lambda.
\end{equation}
Recall  that
\begin{equation}
\label{fixT}
T = \min \big\{  T_0,  \frac{\theta_1 - \theta }{\lambda} \big\} .
\end{equation}
Note that for $t \in [0, T]$, $s(t)$ remains in $[\theta , \theta_1]
\subset ]1- \alpha, \alpha [$.

   We will consider solutions of the Cauchy problem
   \begin{equation}
\label{Cauchyb}
L u = f, \quad  u_{|t = 0} = u_0 , \quad   X u_{| t = 0 }  = u_1
\end{equation}
   with
    \begin{equation}
\label{assu1}
   u \in \cH_{ 1- \theta+ \frac{1}{2} log, \lambda} (T),
     \qquad \D_t u \in \cH_{ - \theta + \frac{1}{2} log, \lambda} (T),
\end{equation}
\begin{equation}
\label{assu2}
    u_0 \in H^{ 1 - \theta} ( \RR^n) ,  \quad      u_1 \in H^{  -
\theta} ( \RR^n),
\end{equation}
\begin{equation}
\label{assu3}
f = f_1 + f_2,    \quad   f_1 \in \cL_{- \theta, \lambda}(T) , \quad
   f_2\in   \cH_{- \theta - \frac{1}{2} log , \lambda}(T),
\end{equation}
Note that if $u$ and $f$ satisfy \eqref{assu1} and \eqref{assu3}, then
\begin{eqnarray}
\label{subassu}
&& u\in L^2([0, T]; H^{1- \theta_1} ) , \quad  \D_t u\in L^2([0, T];
H^{- \theta_1} ),
\\
&& f \in L^1([0, T]; H^{- \theta_2} )
\end{eqnarray}
   for all $ \theta_2 \in ]\theta_1 , \alpha [$, so that the meaning of
the Cauchy condition is clear.

The main step in the proof of Theorem \ref{mainthm} is the following:

\begin{theo}
\label{thm42}
There is a $\lambda_0 \ge 0 $ of the form $\eqref{lda0}$
such that for $\lambda \ge \lambda_0 $
there is a constant $K$ such that:  for all $f$, $u_0$ and $u_1$ satisfying
$\eqref{assu2}$ $\eqref{assu3}$, and all     $u$ satisfying
$\eqref{assu1}$ solution of the Cauchy problem~$\eqref{Cauchyb}$, then
   \begin{equation}
\label{cont}
   u \in \cC_{1- \theta, \lambda}(T), \quad
    \D_t u \in \cC_{- \theta, \lambda}(T)
\end{equation}
       and $u$ satisfies the energy estimate $\eqref{mainest}$.

\end{theo}

This theorem contains two pieces of  information :  first an energy estimate
for smooth $u$,
see Propositions~\ref{prop45cc} and \ref{prop45}.
   By a classical argument, smoothing the coefficients and passing to
the limit,  this estimate allows for  the construction of weak
solutions, see Section 5.2.
The second  piece of information contained in the theorem is a ``weak=strong''
type  result showing that
for data as in the theorem,
    any (weak) solution $u$  satisfying \eqref{assu1} is the limit of
smooth (approximate) solutions, in the norm
    given by the left hand side of the energy estimate,  implying that
$u$ satisfies the additional smoothness \eqref{cont} and   the energy
estimate. This implies uniqueness
of weak solutions.

The idea is to get an energy estimate by integration by parts,
from the analysis of
\begin{equation}
\label{LQ}
2 \re  \langle L u ,  e^{- 2 \gamma t}(1 - \Delta_x)^{ -  s(t)} X  u  \rangle
\end{equation}
where $\langle\cdot, \cdot \rangle $ denotes the $L^2$ scalar product in
$\RR^n$ extended  to the Hermitian symmetric duality  $ H^\sigma
\times H^{- \sigma}$ for all $ \sigma \in \RR$,
and    $  \Delta_x $ denote the Laplace operator on $\RR^d$.
This extends the analysis of \cite{CL} where $X = \D_t$.
   The parameter
$\gamma $ is  chosen  at the end  to absorb classical error terms
(present for Lipschitz coefficients)
while the parameter $\lambda$ which enters in the definition of
$s(t)$, is chosen to absorb
extra error terms coming from the loss of smoothness of the coefficients.

To prove  Theorem~\ref{thm42}, the first idea would be to   mollify
the equation. However, the
lack of smoothness of the coefficients does not allow us to use this
method directly and
we cannot prove that the weak solutions are limits of exact smooth solutions.
Instead, the idea is to write the equation as a system \eqref{syste}
for $(u, v)$ and mollify
this system. This leads to the consideration of  the equations:
\begin{equation}
\label{syste1}
\left\{\begin{aligned}
&  Y^* v + \tilde b_0   v  = \tilde L_2 u - \tilde L_1 u - \tilde d u  + f,
\\
&    Y   u + \tilde c _0 u  =   v / a_0 + g .
\end{aligned}\right.
\end{equation}
   In this form, the commutator of spatial mollifiers with $\D_t$ are trivial,
   and we can prove that weak solutions of  \eqref{syste1} are limits of
   smooth solutions, $(u^\eps, v^\eps)$ with $g^\eps \ne 0$,
   which thus do not correspond to exact solutions $u^\eps$ of \eqref{Cauchyb}.
\medbreak

\noindent\textbf{Notations. } It is important for our purpose to keep track
of the dependence of the various constants on the Log-Lipschitz norms.  In
particular we will use the notations $\delta_0, \delta_1$ of
\eqref{hyp} and $A_{LL}, A_{L^\infty} , B $ of  \eqref{defA}
\eqref{defB}. To simplify the exposition,
we will denote by   $C$, $K_0$  and $K$  constants which may
vary from one line to another,
$C$  denoting universal constants depending only  on the
paradifferential calculus; $K_0$ depending also on
   $A_{L^\infty}/ \delta_0$;  $K$, still independent of the parameters
$(\gamma, \eps)$,   but dependent also on  $\delta_0$, $\delta_1$,
$\theta_0$, $\theta_1$  and the various  norms of the coefficients.

%%%%%%%%%%%%%%%%%%%%%%%%%%

\subsection{Estimating $v$}

First, we give estimates that link $v$ and $\D_t u$.

\begin{lem}
\label{lem43}
   Suppose that $ u $ satisfies $\eqref{assu1}$. Then
       $ v  = X u + c _0 u$ belongs to the space $ \cH_{ - \theta+
\frac{1}{2} log, \lambda} (T)$
$\subset L^2([0, T] ; H^ { - \theta_1})$ and   for almost all $t$,
\begin{equation}
\label{estv1}
\begin{aligned}
\| v(t) \|_{H^{-s(t) + \frac{1}{2}log} } \le &  C
    A_{L^\infty}     \big(  \| u(t) \|_{H^{1 -s(t) + \frac{1}{2}log} } +
\| \D_t u (t) \|_{H^{-s(t) + \frac{1}{2}log} } \big)
\\
   + & C
    (  A_{LL}  + B)     \big(  \| u(t) \|_{H^{1 -s(t) }  } +
\| \D_t  u(t) \|_{H^{-s(t)  } } \big),
\end{aligned}
\end{equation}
\begin{equation}
\label{estv2}
\begin{aligned}
\| \D_t u (t)  \|_{H^{-s(t) + \frac{1}{2}log} }  \le  &
      K_0        \| u(t)  \|_{H^{1 -s(t) + \frac{1}{2}log} }
+
    \frac{C}{\delta_0} \| v  (t) \|_{H^{-s(t) + \frac{1}{2}log} }
\\
& +
K     \big(  \| u(t) \|_{H^{1 -s(t) }  } +
\| v(t) \|_{H^{-s(t)  } } \big).
\end{aligned}
\end{equation}
There are similar estimates in the spaces $H^{s}$ without the $\ml$.

If in addition $Lu = f$ with $f$ satisfying $\eqref{assu3}$,
then
$ \D_t v \in L^1([0, T] ; H^ {-1 - \theta_1})$.
\end{lem}

\begin{proof}  {\bf a) }   First, we note that the multiplication
$ (a, u) \mapsto au$ is continuous from $( L^\infty \cap LL)
([0,T]\times \RR^n)   \times
\cH_{ - \theta + \frac{1}{2} log, \lambda} (T)$ to $\cH_{ - \theta +
\frac{1}{2} log, \lambda} (T)$.
Indeed, the corresponding norm estimate of the product is clear for
smooth $u$, from
\eqref{estprod} integrated in time. The claim follows by density.
In particular, this shows that  $ a_0 \D_t u$ and the $ a_j \D_{x_j}
u $  belong to $ \cH_{ - \theta+ \frac{1}{2} log, \lambda} (T)$.
Similarly, the estimate
\begin{equation}
\| b u (t) \|_{H^{-s(t) + \frac{1}{2}log} } \le  C \| b u (t)
\|_{H^{1 -s(t) }  }
   \le C  \| b \|_{C^\alpha}   \| u(t) \|_{H^{1 -s(t) }  }
\end{equation}
   implies that $ c_0 u  \in  \cH_{ - \theta+ \frac{1}{2} log, \lambda} (T)$.
Therefore  $ v \in \cH_{ - \theta+ \frac{1}{2} log, \lambda} (T)$ and
the  estimate
\eqref{estv1} holds. The proof of \eqref{estv2} is similar, noting that
$$
\D_t u = \frac{1}{a_0} v  - \sum_{j=1}^d  \frac{a_j}{a_0} \D_{x_j } u
- \frac{c_0}{a_0} u .
$$

{\bf b) } As in the proof of Lemma~\ref{lem42}, we see that the
equation implies that
\begin{equation*}
\D_t  v   = f   -  \sum_{j=1}^d  \D_{x_j} ( \tilde a_j v) - \tilde
b_0 v   + \tilde L_2 u - \tilde L_1 u -  \tilde d u .
\end{equation*}
The conservative form of  $ \tilde L_2$ and the multiplicative
properties  above show that
$$
\D_{x_j} ( \tilde a_j v) , \ \tilde L_2 u \in \cH_{ - \theta -1 +
\frac{1}{2} log, \lambda} (T)
\subset L^2([0, T] ; H^ {-1 - \theta_1}).
$$
Similarly, $\tilde L_1  u $ and $ \tilde b_0 v$   belong  to
$  \cH_{ - \theta  + \frac{1}{2} log, \lambda} (T) $ and thus to
   $ L^2([0, T] ; H^ { - \theta_1}) $.
The last term $ \tilde  u $ is in $ L^2 $.
Therefore, $\D_t v -  f  \in L^2([0, T] ; H^ {-1 - \theta_1})$. Since
$f \in  L^1([0, T] ; H^ { - \theta_2 }) $  for $\theta_2 \in
]\theta_1, \alpha [$, the lemma follows.
\end{proof}

Next, we give a-priori estimates in the space
$ \cH_{ - \theta+ \frac{1}{2} log, \lambda} (T) \cap \cC_{- \theta,
\lambda}(T) $  for
smooth solutions of
\begin{equation}
( Y^* + \tilde c_0) v =    \varphi  , \quad  v_{| t =0}  =v_0.
\end{equation}
We define the operators
   \begin{equation}
   \label{defQ}
( Q v) (t)   = (1 - \Delta_x)^{- s(t)/2}  v(t) , \quad
   ( Q_\gamma  v) (t) =  e^{ - \gamma t } (Q v)(t).
   \end{equation}

\begin{prop}
\label{prop45cc}
Suppose that   $v \in L^2([0, T]; H^1)$   and
$\D_t v \in L^1 ([0, T]; L^2)$.
   Then  the functions
       $v_{\gamma}  (t) :=   Q_\gamma  v $ belong to $ C^{0}([0, T],
L^2) $    and satisfy
       \begin{equation}
\label{410a}
       \begin{aligned}
       \| v_\gamma  (t)  \|^2_{L^2}&  +  2  \int_{0}^t \|
       (\gamma +
       \lambda \Lambda )^{1/2} v_{\gamma} (t') \|^2_{L^2} dt'
       \\
         \le &  2  \int_0^t
    \langle (Y^{*} + \tilde c_0) v (t')   ,   Q_{\gamma}^{2}(t')   v
(t')    \rangle   \, dt'  +
     \| v_{\gamma, \eps}(0) \|^2_{L^2}
     \\
    & +   \int_{0}^t F(t') dt'
\end{aligned}
\end{equation}
with
\begin{equation}
\label{errv}
    F (t')  \le
    K_0   \frac{ A_{LL} }{\delta_0}   \| e^{-\gamma {t'}}
\Lambda^{1/2}  v(t') \|^{2 }_{H^{-s(t')}}   + K   \| v(t' ) \|^2
_{H^{ - s(t')}}    .
\end{equation}
\end{prop}

\begin{proof}

{\bf a)}  Since  $ v \in   L^2 ([0,T] ; H^{1}) $
and $ \D_t v \in L^1([0, T] ; L^2)$, we have
\begin{equation}
\label{comQDt}
\D_t Q_{\gamma}  v  =    Q_{\gamma} \D_t  v    -
   ( \gamma + \lambda \Lambda )  Q_{\gamma}   v  \in L^1([0,T]; L^2 )
\end{equation}
as immediately seen using the spatial Fourier transform.
Moreover,
   $ v_\gamma = Q_{\gamma}v  \in  C^0 ([0,T]; L^2 ) $
and satisfies the following identity
$$
   \| v_{\gamma}(t)  \|^2_{L^2} -  \| v_{\gamma}(0)  \|^2_{L^2} =
2 { \rm Re}  \int_0^t \langle \D_t Q_{\gamma}     v   ,   Q_{\gamma}
v   \rangle  \, dt' .
$$
Thus,
\begin{equation}
\label{qqq}
\begin{aligned}
2&  {\rm Re}  \int_0^t   \langle \D_t  v   ,   Q_{\gamma}^{2}  v
\rangle  \, dt' =
2 { \rm Re} \int_0^t   \langle Q_{\gamma}  \D_t  v   ,   Q_{\gamma}
v   \rangle  \, dt'
\\ & =
   \| v_{\gamma}(t)  \|^2_{L^2} -  \| v_{\gamma}(0)  \|^2_{L^2}
        +  2
    \int_{0}^t \|    (\gamma +
       \lambda \Lambda )^{1/2} v_{\gamma} (t') \|^2_{L^2} dt'
       \end{aligned}
\end{equation}

\medbreak
{\bf b)}  Next we consider the terms $ \D_{x_j} (\tilde a_j   v )$.
We note that they belong to $L^2([0, T]; H^{- \sigma} )$ for all $\sigma  > 0$.
In particular, since $s (t) \ge \theta > 0$, we note that the pairing
$$
\langle \D_{x_j} (\tilde a_j   v ) , Q_\gamma ^2 v \ra
$$
is well defined.
We give an estimate for
$$
   2 { \rm Re}  \int_0^t \     \langle \D_{x_j} (\tilde a_j   v )   ,
Q_{\gamma}^{2}  v   \rangle  dt' ,
$$ using the decomposition
$$
\tilde a_j   v  =  T_{\tilde a_j   } v  + R_{\tilde a_j  } v .
$$
By Proposition~\ref{proppara} it follows
$$
\| R_{\tilde a_j }   v(t)  \| _{H^{ 1 - s(t) -\frac{1}{2} log}}   \le
C \| \tilde a_j\|_{LL}
\|     v(t) \| _{H^{  - s(t) + \frac{1}{2} log}}
$$
since $ s(t) \in [ \theta, \theta_1] \subset]0, 1[$.  Moreover,
$$
\|    Q^2_{\gamma}  v(t) \| _{H^{   s(t) + \frac{1}{2} log}}   \le  C
e^{ - 2 \gamma t}  \|     v(t) \| _{H^{  - s(t) + \frac{1}{2} log}} .
$$
Thus
\begin{equation*}
\begin{aligned}
   |  \langle \D_{x_j} R_{\tilde a_j }   v  (t)  ,    Q_{\gamma,
\eps}^{2}  v (t)   \rangle   |
&  \le
\| R_{\tilde a_j }   v(t)  \| _{H^{ 1 - s(t) -\frac{1}{2} log}}   \
\|    Q^2_{\gamma, \eps}  v(t) \| _{H^{   s(t) + \frac{1}{2} log}}
\\
   & \le
C  \| \tilde a_j\|_{LL}  e^{ - 2 \gamma t}  \|     v(t) \|^2  _{H^{
- s(t) + \frac{1}{2} log}}.
\end{aligned}
\end{equation*}
It  remains to consider
$$
\begin{aligned}
{\rm Re}  \langle \D_{x_j} T_{\tilde a_j}   v    ,   Q_{\gamma}^{2}
v   \rangle  &  =
{\rm Re}  \langle Q_{\gamma} \D_{x_j} T_{\tilde a_j}   v    ,
Q_{\gamma}  v   \rangle
\\
& = {\rm Re}  \langle  \D_{x_j} T_{\tilde a_j}   Q_{\gamma}  v    ,
Q_{\gamma}  v   \rangle
+ {\rm Re}  \langle \D_{x_j}  [ Q_{\gamma} , T_{\tilde a_j} ]   v
,   Q_{\gamma}  v   \rangle .
\end{aligned}
$$
Note that these computations make sense because
$ v(t)  \in  H^ {1} $ and   all the pairings are well defined.
   Proposition~\ref{propB3b} implies that
$$
\|  \langle   \D_{x_j}  [ Q_{\gamma}, T_{\tilde a_j} ]   v (t) \|_{0-\mez log}
\le C  e^{ - \gamma t} \| \tilde a_j \|_{LL} \|   v (t) \|_{- s(t) + \ml }
$$
and therefore
\begin{equation}
\label{415a}
   |  \langle   \D_{x_j}  [ Q_{\gamma}, T_{\tilde a_j} ]   v (t)   ,
Q_{\gamma}  v  (t) \rangle   |
\le
C  \| \tilde a_j\|_{LL}  e^{ - 2 \gamma t}  \|     v(t) \|^2  _{H^{
- s(t) + \frac{1}{2} log}} .
\end{equation}
Next, for $ v_{\gamma } (t) \in H^{2 - \theta_1} $, we have
$$
\begin{aligned}
2 {\rm Re}  \langle  \D_{x_j} T_{\tilde a_j}   v_{\gamma}       ,   &
v_{\gamma}      \rangle
= {\rm Re}  \langle  (\D_{x_j} T_{\tilde a_j}   -  T^* _{\tilde a_j}
\D_{x_j})  v_{\gamma}       ,   v_{\gamma}      \rangle
\\
& =
{\rm Re}  \langle  ( T_{\tilde a_j}   -  T^* _{\tilde a_j} ) \D_{x_j}
v_{\gamma}       ,   v_{\gamma}      \rangle
+ {\rm Re}  \langle [  \D_{x_j} , T_{\tilde a_j} ]    \D_{x_j}
v_{\gamma,}       ,   v_{\gamma}      \rangle   .
\end{aligned}
$$
   Using Propositions~\ref{prop37} and \ref{propB3b}, one can bound
both terms by the right hand side of
\eqref{415a}.
Adding up, we have proved that
\begin{equation*}
\label{qqq2}
\begin{aligned}
\big|  2 {\rm Re}  \int_0^t        \langle \D_{x_j} (\tilde a_j   v )
,   Q_{\gamma}^{2}  v   \rangle   \, dt'  \big| \le
      C \| \tilde a_j \|_{LL}
    \int_{0}^t \|   e^{- \gamma {t'}}
        \Lambda ^{1/2} v  (t') \|^2_{H^{- s(t)}} dt'.
       \end{aligned}
\end{equation*}

   \medbreak
   {\bf c)} The zero-th order term is clearly a remainder, and  the
multiplicative properties
   imply that
   $$
   | \la \tilde c_0 v (t)  , Q_\gamma^2 v  (t) \ra \le
K     \| v (t) \|^2 _{ H^{ - s(t)}  } .
   $$

   \medbreak
{\bf  d)}  We note that
$$
\begin{aligned}
\| a_j/ a_0 \|_{LL} &  \le \| a_j\|_{LL} \| \| 1/  a_0 \|_ {L^\infty}
+ \| a_j\|_{L^\infty}  \| \| 1/  a_0 \|_  {LL}
\\
& \le  \frac{ A_{LL}}{\delta_0}   +  \frac{ A_{L^\infty} A_{LL}}{\delta^2 _0}
\le 2 \frac{ A_{L^\infty} A_{LL}}{\delta^2 _0}  ,
\end{aligned}
$$
since  $ \delta_0 \le a_0 \le A_{L^\infty}$.
Using identity \eqref{qqq} and  the estimates of parts b)  and
c), implies \eqref{410a}  and so the proof of the Lemma is complete.
\end{proof}

%%%%%%%%%%%%%%%%%%%%%%%%%

\subsection{Estimating $\nabla_x u$}

We now get estimates of  $\nabla_x u$ from the analysis of
\begin{equation}
\label{417cc}
-  2  \re \la \widetilde L_2 u, Q_{\gamma}^{2}  X u    \rangle =
    - \sum_{j, k= 1}^n  2 \re \la \D_{x_j} (  \tilde a_{j, k} \D_{x_k}
u )    ,   Q_{\gamma}^{2} Xu    \rangle
\end{equation}

   \begin{prop}
   \label{prop45}
   Suppose that $u \in L^2([0, T]; H^2)$ with $\D_t u \in L^2([0, T]; H^1)$.
   Then
       $u_{\gamma} := Q_{\gamma} u \in C^{0}([0, T], H^{1})  $    and
       \begin{equation}
\label{418c}
       \begin{aligned}
    \mez   &  \delta_0 \delta_1     \| \na_x    u_{\gamma}(t)
\|^2_{L^2}   + \int_{0}^t
   \delta_0 \delta_1   \|(\gamma +
       \lambda \Lambda )^{1/2} \na_x u_{\gamma} (t') \|^2_{L^2} dt'
       \\
         \le &   -  2 {\rm Re}
   \int_0^t   \la \widetilde L_2 u, Q_{\gamma}^{2}  v   \rangle
\, dt'    +
   C A^2_{L^\infty}   \| \na_x u_{\gamma}(0) \|^2_{L^2}
+   \int_{0}^t  E( t')  dt',
\end{aligned}
\end{equation}
where
\begin{equation}
\label{err}
\begin{aligned}
   | E(t) |& \le
   \\
    &  K_0
      A_{LL} A_{L^\infty}  e^{-2\gamma t}   \big(   \| u(t) \|^{2
}_{H^{1 -s(t)+\ml }}   +  \frac{1}{\delta_0^2}
   \| Xu (t)  \|^{2 }_{H^{ -s(t)+\ml }} \big)
       \\
     &  +  K  e^{-2\gamma t}    \big(   \| u(t) \|^{2 }_{H^{1 -s(t)}}   +
   \| Xu (t)  \|^{2 }_{H^{ -s(t) }} \big)  .
\end{aligned}
\end{equation}

   \end{prop}

To simplify the exposition, we note here  that all the dualities $
\langle \cdot, \cdot \rangle $ written
below make sense, thanks to the  smoothness assumption on $u$.
This will not be repeated  at each step.
Moreover, in the proof below, we assume that $u$ itself is smooth (in time).

   \begin{proof} {\bf a)}  We first perform several reductions.
Using $iii)$ of Proposition~\ref{proppara}, one shows that
$$
    \la \D_{x_j} (  \tilde a_{j, k} \D_{x_k} u )    ,   Q_{\gamma}^{2}
Xu    \ra =
\la  \D_{x_j}  ( T_{ \tilde a_{j, k} } \D_{x_k} u   )   ,
Q_{\gamma}^{2}  Xu    \ra  +  E_1
$$
with
\begin{equation}
\label{err1}
|  E_1 (t) |  \le C  \| \tilde a_{j, k} \|_{LL}  \| \D_{x_k} u(t)
\|_{H^{ - s(t)+ \frac{1}{2}log}}
\| Q^2_{\gamma} Xu (t) \|_{H^{  s(t)+ \frac{1}{2}log}} .
\end{equation}
Since $   \| \tilde a_{j, k} \|_{LL}  \le K_0 A_{LL} \le K_0 A_{LL}
A_{L^\infty} / \delta_0 $,
$E_1 $ satisfies \eqref{err}.
   Similarly,
$$
\begin{aligned}
\la  \D_{x_j}  ( T_{ \tilde a_{j, k} } \D_{x_k} u   )   ,
Q_{\gamma}^{2}  Xu   \ra
& =
\la   \D_{x_j}   Q_{\gamma, \eps}  T_{ \tilde a_{j, k} } \D_{x_k} u
,     Q_{\gamma}  Xu   \ra  +
\\
& =    \la   \D_{x_j}     T_{ \tilde a_{j, k} } \D_{x_k} Q_{\gamma} u      ,
      Q_{\gamma}  X u   \ra  +
E_2
\end{aligned}
$$
where $E_2$ also satisfies \eqref{err1}, and hence \eqref{err}.

\medbreak
{\bf b)}  Next we write
\begin{equation*}
X u   = T_{a_0} \D_t u  + \sum T_{a_j} \D_{x_j} u  + r
\end{equation*}
and
   \begin{equation*}
\begin{aligned}
\| r(t) \|_{ H^{1-s(t) -\ml} } \le C A_{LL} \big( \| u(t) \|_{
H^{1-s(t) + \ml}} + \| \D_t u (t) \|_{ H^{-s(t) + \ml}}
\big)
\\
+  C  B \| u (t) \|_{H^{1 - s(t)}}.
\end{aligned}
\end{equation*}
   Therefore,  $r$ contributes to an error term
$ E_3 =  \la  \D_{x_j}     T_{ \tilde a_{j, k} } \D_{x_k} Q_{\gamma}
u      ,     Q_{\gamma}  r \ra$
such that
\begin{equation*}
| E_3(t) | \le   e^{ -2 \gamma t}  K_0 A_{L^\infty}
\| u (t)  \|_{H^{1-s(t)+\ml}}  \| r(t)  \|_{H^{1-s(t)-\ml}}  .
\end{equation*}
Using \eqref{estv2} in the estimate of $r$, we see that
\begin{equation*}
\begin{aligned}
| E_3(t) | \le   e^{ -2 \gamma t}    K_0 A_{L^\infty} A_{LL}
\| u (t) & \|_{H^{1-s(t)+\ml}}
\\
   \Big( \| u (t)  \|_{H^{1-s(t)+\ml}}   & +
\frac{1}{\delta_0} \| Xu  (t)  \|_{H^{-s(t)+\ml}}
\\ & + K \| u (t) \|_{H^{1 - s(t)}}
+  K \| X u (t) \|_{H^{ - s(t)}} \Big)
\end{aligned}
\end{equation*}
and hence satisfies \eqref{err}.

\medbreak
{\bf c)}  Consider now the term
\begin{equation*}
\begin{aligned}
   \la  \D_{x_j}     T_{ \tilde a_{j, k} } \D_{x_k} Q_{\gamma} u      ,  \,
&   Q_{\gamma}  T_{a_0} \D_t u \ra  =
-   \la  T_{ \tilde a_{j, k} } \D_{x_k} Q_{\gamma} u      ,
\D_{x_j}   Q_{\gamma}  T_{a_0} \D_t u \ra
\\
& = -   \la  T_{ \tilde a_{j, k} } \D_{x_k} Q_{\gamma} u      ,
T_{a_0}    \D_{x_j}   Q_{\gamma}  \D_t u \ra  + E_4
\\
& = -  \la  (T_{a_0}) ^*T_{ \tilde a_{j, k} } \D_{x_k} Q_{\gamma} u
,       \D_{x_j}   Q_{\gamma}  \D_t u \ra  + E_4
\\
& = -  \la  T_{a_0} T_{ \tilde a_{j, k} } \D_{x_k} Q_{\gamma} u
,       \D_{x_j}   Q_{\gamma}  \D_t u \ra  + E_4 + E_5
\\
& = -  \la  T_{a_0  \tilde a_{j, k} } \D_{x_k} Q_{\gamma} u      ,
\D_{x_j}   Q_{\gamma}  \D_t u \ra  + E_4 + E_5 + E_6
\end{aligned}
\end{equation*}
where $E_4$, $E_5$ and $E_6$  are estimated by
Proposition~\ref{propB3b}, \ref{prop37} and \ref{prop38}
respectively. They all satisfy
\begin{equation*}
| E_k (t) | \le C    e^{ - 2 \gamma t}  A
\| u (t)  \|_{H^{1-s(t)+\ml}}  \| \D_t  u (t)  \|_{H^{-s(t)+\ml}}  .
\end{equation*}
with
$ A = \| \tilde a_{j, k}\| _{LL} \| a_0 \|_{L^\infty} + \| \tilde
a_{j, k}\|_{L^\infty} \| a_0 \|_{LL}
\le K_0 A_{L^\infty} A_{LL}$.
Again using \eqref{estv2} to replace $\D_t u$ by $Xu $, one shows
that these errors satisfy \eqref{err}.

Similarly
\begin{equation*}
\begin{aligned}
   \la  \D_{x_j}     T_{ \tilde a_{j, k} } \D_{x_k} Q_{\gamma} u ,   & \,
Q_{\gamma}  T_{a_l } \D_{x_l}  u \ra
\\
& = -  \la  T_{a_l  \tilde a_{j, k} } \D_{x_k} Q_{\gamma} u      ,
\D_{x_l}    \D_{x_j}   Q_{\gamma}    u \ra  + E_7
\end{aligned}
\end{equation*}
where $E_7$ satisfies
\begin{equation}
\label{err4}
| E_7(t) | \le C    e^{ - 2 \gamma t}  K_0 A_{L^\infty} A_{LL}
\| u (t)  \|^2 _{H^{1-s(t)+\ml}}
\end{equation}
   thus \eqref{err}.

\medbreak
{\bf d)} Introduce the notation
\begin{equation}\label{422c}
w_j = \D_{x_j}   Q_{\gamma}    u.
\end{equation}
Because $\tilde a_{j,k} = \tilde a_{k, j}$, we have
\begin{equation*}
\begin{aligned}
    \re \la  T_{a_l  \tilde a_{j, k} } w_k , \D_{x_l} w_j \ra & +
     \re \la  T_{a_l  \tilde a_{k, j} } w_j  , \D_{x_l} w_k \ra
\\
& =   \re \la ( (T_{a_l  \tilde a_{j, k} })^* \D_{x_l} -  \D_{x_l}
T_{a_l  \tilde a_{j, k} }  )w_k ,  w_j \ra := E_8
\end{aligned}
\end{equation*}
Using Propositions~\ref{prop37}  and \ref{propB3b}, one shows that 
$E_8$ satisfies
$$
| E_8(t) |  \le C \| a_l \tilde a_{j, k} \| _{LL}   \| w_j (t)
\|_{H^{0+\ml}}  \,   \| w_k (t) \|_{H^{0+\ml}}
$$
and therefore $E_8$ also satisfies \eqref{err4} thus \eqref{err}.

\medbreak
{\bf e)}  It remains to consider the sum
\begin{equation}
S :=   \re  \sum_{j, k = 1}^n  \la  T_{b_{j, k} } \D_{x_k}
Q_{\gamma} u      ,       \D_{x_j}   Q_{\gamma}  \D_t u \ra
\end{equation}
with $b_{j, k} = a_0  \tilde a_{j, k}= a_0 a_{, k} + a_j a_k$.
   By the strict  hyperbolicity assumption \eqref{hyp},
   it follows for all
$\xi \in \RR^n$
$$
    \sum_{j, k = 1}^n b_{j, k} (t, x)   \xi_j \xi_k  \ge
\delta_0\delta_1  | \xi |^2 .
$$
Therefore, we can use  Corollary~\ref{cor316b}.
Since $ \| b_{j, k} \|_{LL}  \le 2 A_{L^\infty} A_{LL}$,   there
exists an integer $\nu$, with
\begin{equation}
\label{425c}
\frac{2^\nu }{\nu}  \approx  \frac{A_{L^\infty} A_{LL}}{   \delta},
\end{equation}
such that for all $ t \in [0, T_0]$ and
$(w_1, \ldots, w_n)$ in $L^2(\RR^n)$,  the following estimate is satisfied
\begin{equation}
\label{Gard}
   \re  \sum_{j, k = 1}^n  \la  P^\nu_{b_{j, k}(t)  } w_k  , w_j \ra
\ge \frac{ \delta_0\delta_1 }{2}  \| w \|^2 _{L^2}
\end{equation}
 From now on we fix such a $\nu$ and use the notation $P_{b} $ in
place of $P^\nu_{b} $.

Using Lemma~\ref{lem312a}  and Proposition~\ref{prop314c}, we see that
$$
\begin{aligned}
   \|   \D_{x_j}     T_{b_{j, k} } w_k  -   \D_{x_j}   \widetilde
P_{b_{j, k} } w_k \|_{H^{0-\ml}}
    \le C\| b_{j, k} \|_{LL} \Big( \| w_k \|_{H^{0+\ml}} +
K  \| w_k \|_{L^2} \Big)
\end{aligned}
$$
Therefore
\begin{equation*}
S  = \re \sum_{j, k = 1}^n  \la  \widetilde P_{b_{j, k} }
\D_{x_k} Q_{\gamma} u      ,       \D_{x_j}   Q_{\gamma}  \D_t u \ra
+ E_9
\end{equation*}
where
\begin{equation*}
\begin{aligned}
| E_9(t)   | \le C    e^{ - 2 \gamma t}  & \| b_{j, k} \|_{LL}
\| u (t)  \|_{H^{1-s(t)+\ml}}  \| \D_t  u (t)  \|_{H^{-s(t)+\ml}}
\\
& +       e^{ - 2 \gamma t} \nu  K
\| u (t)  \|_{H^{1-s(t)}}  \| \D_t  u (t)  \|_{H^{-s(t)+\ml}}  .
\end{aligned}
\end{equation*}
Using \eqref{estv2}, implies that $E_9$ satisfies \eqref{err}.

Next, we use Proposition \ref{prop314c}  to replace
$ \D_{x_j}  \widetilde P_{b_{j, k} }$  by
$ \mez   \D_{x_j}  ( \widetilde P_{b_{j, k} }+  (\widetilde  P_{b_{j,
k} })^* ) $
at the cost of an error $E_{10}$  similar to $E_9$.

At this  stage, we commute $Q_{\gamma}$ and $  \D_t $ as in
\eqref{comQDt}. Using the notation \eqref{422c},   yields
\begin{equation}
\begin{aligned}
2 S    =  &     \sum_{j, k = 1}^n
\re  \la  ( \widetilde P_{b_{j, k} 
} +  (\widetilde  P_{b_{j, k} })^*
) w_k      ,        \D_t w_j \ra
\\
& + \gamma     \sum_{j, k = 1}^n
\re  \la  ( \widetilde P_{b_{j, k} } +  (\widetilde  P_{b_{j, k} })^*
) w_k      ,        w_j \ra
\\
& + \lambda    \sum_{j, k = 1}^n
\re  \la  ( \widetilde P_{b_{j, k} } +  (\widetilde  P_{b_{j, k} })^*
) w_k      ,    \Lambda     w_j \ra
+ 2 E_9 + 2 E_{10} .
\end{aligned}
\end{equation}
   We denote by $S^1$, $S^{2}$ and $S^{3} $  the   sums  on the right hand side.
\medbreak

{\bf f)}  The symmetry $ b_{j, k} = b_{k,j}$ implies the identity
\begin{equation*}
S^ 1  =    \frac{d}{dt}     \sum_{j, k = 1}^n \re
    \la   \widetilde P_{b_{j, k} }   w_k      ,      w_j \ra         +   E_{11}
\end{equation*}
where
\begin{equation*}
E_{11}  =       \sum_{j, k = 1}^n \re
    \la  [  \widetilde   P_{b_{j, k} }  , \D_t ]  w_k      ,      w_j \ra
\end{equation*}
is estimated  using  Proposition~\ref{prop314c}:
\begin{equation*}\label{err6}
\begin{aligned}
|  E_{11} (t)   |     \le  C     \| b_{j, k} \|_{LL}  &\big(
\| w (t)  \|_{H^{0+\ml}}     + \nu \| w \|_{L^2} \big) \| w (t)  \|_{H^{0+\ml}}
\\
     \le   C    e^{ - 2 \gamma t}    \| b_{j, k} \|_{LL} &
    \| u (t)  \|_{H^{1- s(t) +\ml}}
\\
& \big(   \| u (t)   \|_{H^{1- s(t) +\ml}}
+ \nu \| u(t)  \|_{H^{1-s(t)} } \big)
\end{aligned}
\end{equation*}
and therefore satisfies \eqref{err}.
Moreover,
\begin{equation*}
\re  \la   \widetilde P_{b_{j, k} }   w_k      ,    \Lambda     w_j \ra  =
\re  \la  \widetilde P_{b_{j, k} }  \Lambda^\mez w_k      ,
\Lambda^\mez      w_j \ra
+ \re  \la   \Lambda^\mez   [  \Lambda^\mez , \widetilde P_{b_{j, k}
},   ] w_k      ,       w_j \ra
\end{equation*}
\begin{equation*}
\re  \la  ( \widetilde P_{b_{j, k} } )^*   w_k      ,    \Lambda     w_j \ra  =
\re  \la  \Lambda^\mez w_k      ,    \widetilde P_{b_{j, k} }
\Lambda^\mez      w_j \ra
+ \re  \la   w_k      ,     [ \widetilde P_{b_{j, k} },
\Lambda^\mez ]   \Lambda^\mez   w_j \ra  .
\end{equation*}
We use Proposition~\ref{prop317} to estimate the commutators and
\begin{equation*}
S^3 = 2   \sum_{j, k = 1}^n \re  \la  \widetilde P_{b_{j, k} }
\Lambda^\mez w_k      ,    \Lambda^\mez      w_j \ra     + E_{12}
\end{equation*}
where
\begin{equation*}
\label{err7}
|  E_{12} (t) |  \le  K   \| w(t)  \|^2_{L^2}
   \le K \| u(t)  \|^2_{H^{1- s(t)}}.
\end{equation*}

Summing up, we have shown that up to an error which satisfies
\eqref{err},  the quantity
\eqref{417cc}
under consideration is equal to
\begin{equation}\label{433c}
\begin{aligned}
   \frac{d}{dt}     \sum_{j, k = 1}^n \re
    \la   \widetilde P_{b_{j, k} }   w_k      ,      w_j \ra
& +   \gamma      \sum_{j, k = 1}^n
   2  \re  \la   \widetilde P_{b_{j, k} }   w_k      ,        w_j \ra
\\
& +   \lambda   \sum_{j, k = 1}^n   2  \re  \la  \widetilde
P_{b_{j, k} }  \Lambda^\mez w_k      ,    \Lambda^\mez      w_j \ra.
\end{aligned}
\end{equation}
By \eqref{Gard}, the last two sums are larger than or equal to
$ \delta_0 \delta _1 \| w (t) \|^2_{L^2}$  and $ \delta_0 \delta_1
\| w (t) \|^2_{H^{0+\ml}}$, respectively.
Similarly, integrating the first term between $0$ and $t$ and using
\eqref{Gard} gives
control of $ \frac{\delta_0 \delta_1 }{2} \| w ( t ) \|_{L^2}$,
finishing the proof of \eqref{418c}.
   \end{proof}

%%%%%%%%%%%%%%%%%%%%%%%%%%%%%%%%

\subsection{A-priori  estimates for the solutions of \eqref{syste1}}

   The proof of Theorem~\ref{thm42} is based on a-priori estimates for smooth
   solutions of the system~\eqref{syste1}.

\begin{theo}
\label{theo42b}
There are   $\lambda_0 \ge 0 $ of the form $\eqref{lda0}$   and $\gamma_0$
such that for $\lambda \ge \lambda_0 $   and $\gamma \ge \gamma_0$
the following
is true:

for all  $u \in L^2([0, T]; H^2)$   and  $v \in L^2([0, T]; H^1)$
with $\D_t u \in L^2([0, T]; H^1)$ and
$\D_t v \in L^1 ([0, T]; L^2)$ and for all parameters
$ \lambda $,  $\gamma$ and  all $t \le T$, the following 
holds:
   \begin{equation}
\label{mainest22}
\begin{aligned}
&  \sup_{   0 \le t' \le t}   e^{ -2 \gamma t'}
    \Big(  \mez      \delta_0 \delta_1   \| u (t')  \|^2_{ H^{1 -s( t')}} +
     \| v  (t')  \|^2_{ H^{-s(t')}} \Big)  \\
   &  + \delta_0 \delta _1 \int_0^t
   e^{ - 2 \gamma t'}     (\lambda    \| u (t' ) \|^2 _{  H^{1 - s( t')
+ \frac{1}{2} log }} +
    \gamma      \| u (t' ) \|^2 _{  H^{1 - s( t')  }} ) dt'
    \\
    & +   \int_0^t
   e^{ - 2 \gamma t'}     (\lambda
    \| v  (t' ) \|^2 _{  H^{ - s( t') + \frac{1}{2} log }}
    +   \gamma      \| v (t' ) \|^2 _{  H^{ - s( t')  }} ) dt'
    \\
    & \le  C    A_{L^\infty}^2    \| u (0)  \|^2_{ H^{1 - \theta }} +
       \|  v(0)   \|^2_{ H^{-\theta }}
      +    2 {\rm Re}
   \int_0^t    \la f  , Q_{\gamma}^{2}  v   \rangle \, dt'
    \\
    &   +  K  \int_0^t  e^{-2 \gamma t'}  \| g (t') \| _{1 - s(t) -  \ml }
      \| u (t') \| _{1 - s(t) +  \ml } dt'  ,
    \end{aligned}
\end{equation}
with  $f = Y^* v + \tilde b_0 v - \tilde L_2 u + \tilde L_1 u + \tilde d u
\in L^1([0, T]; H^{\alpha' -1} ) $,
$g = Y u + \tilde c_0 u - v/ a_0 \in L^2([0, T]; H^{\alpha' } ) $  for
all  $\alpha' < \alpha$.

      \end{theo}

\begin{proof}
   We compute the integral over $[0, t] $ of $  \re \la f , Q_\gamma^2
v \ra $.  Proposition~\ref{prop45cc}
takes care of the  first term $ 2 \re \la Y^* v + \tilde b_0 v,
Q_\gamma^2 v \ra $.
We split the   second term into three pieces
$$
    \la   \tilde L_2 u , Q_\gamma^2 v \ra =
    \la   \tilde L_2 u , Q_\gamma^2 X u  \ra    -
    \la   \tilde L_2 u , Q_\gamma^2  (a_0 g )  \ra +
   \la   \tilde L_2 u , Q_\gamma^2 (c_0 u)  \ra
$$
   and use Proposition~\ref{prop45} for the first piece.
   The multiplicative properties imply that
   $$
   \begin{aligned}
     |  \la   \tilde L_2 u (t)  , Q_\gamma^2  (a_0 g ) (t)  \ra | \le
     & K   \| g (t) \| _{1 - s(t) -  \ml }  \| \tilde L_2 u (t ) \|_{-
1 - s(t) + \ml }
     \\
     \le   &  K   \| g (t) \| _{1 - s(t) -  \ml }  \|    u (t ) \|_{ 1
- s(t) + \ml },
\end{aligned}
   $$
   and
   $$
   \begin{aligned}
     |  \la   \tilde L_2 u (t)  , Q_\gamma^2  (c_0 u ) (t)  \ra | \le
     & K   \| u (t) \| _{1 - s(t)  }  \| \tilde L_2 u (t ) \|_{- 1 - s(t) }
     \\
     \le   &  K     \|   u (t ) \|^2_{ 1 - s(t)  }.
\end{aligned}
   $$
   Next, using the multiplicative properties stated in
Corollary~\ref{cor36}  for the products
    $ \tilde b_j \D_{x_j} u$ and   $ \D_{x_j} ( \tilde c_j  u) $,
and the embedding $ L^2 \subset H^{ - s} $ for the term $\tilde d u$,
we see that
   \begin{equation*}
\label{432d}
\| ( \widetilde L_1 u + \tilde d  u)(t)  \|_{H^{-s(t)}}  \le K  \|
u (t)  \|_{H^{1-s(t)}}  .
\end{equation*}
Thus
   $$
    \begin{aligned}
     |  \la  ( \tilde L_1 + \tilde d)u (t)  , Q_\gamma^2  v (t)  \ra | \le
     &  K    \| u (t) \| _{1 - s(t)  }  \|  v  (t ) \|_{  - s(t) }
     \\
     \le   &  K    \big(  \|   u (t ) \|^2_{ 1 - s(t)  } + \|  v  (t )
\|^2_{  - s(t) }\big).
\end{aligned}
   $$
    Proposition~\ref{prop45} gives an estimate of $\na_x u$. We also
need an estimate for $u$.
The identity  \eqref{qqq} applied to $u$ yields
\begin{equation*}
\begin{aligned}
e^{ - 2 \gamma t}  \| u  (t) \|^2_{H^{- s(t)} }   +
    \int_0^t
   e^{ - 2 \gamma t'}     (\lambda    \|   u (t' ) \|^2 _{  H^{ - s(
t') + \frac{1}{2} log }} +
    \gamma      \|   u (t' ) \|^2 _{  H^{ - s( t')  }} ) dt'
   \\
    \qquad =   \| u_{\gamma}  (0) \|^2_{H^{-s(0)}}   +
2  \re  \int_0^t  \la  \D_t u , Q_\gamma^2 u \ra  dt' .
\end{aligned}
\end{equation*}
Next, we use the inequality
$$
   | \la  \D_t u , Q_\gamma^2 u \ra |  \le C  \big(  \| u (t) \|^2_{H^{
1 - s(t)}}  +
   \| \D_t u (t) \|^2_{H^{-1  - s(t)}}. 
    \big)
$$
In addition, we note that  the second equation in \eqref{syste1} yields
$$
\| \D_t u (t) \|_{H^{-1  - s(t)}}  \le  K \big( \| v  (t) \|^2_{H^{
- s(t)}}  +
\|   u (t) \|^2_{H^{  - s(t)}} \big)  + \| g (t) \|^2_{H^{-1  - s(t)}}.
$$

We add the various estimates and use Propositions~\ref{prop45cc} and
   \ref{prop45} to obtain a final estimate.
    On the left hand side we have
\begin{equation}
\label{LH1}
      \sup_{   0 \le t' \le t}   e^{ -2 \gamma t'}
    \Big(  \mez      \delta_0 \delta_1   \|  u (t')  \|^2_{ H^{ 1-s( t')}} +
     \| v  (t')  \|^2_{ H^{-s(t')}} \Big)
      \end{equation}
\begin{equation}
\label{LH2}
+     \gamma \int_0^t   e^{ -2 \gamma t'}
\big(  \delta_0 \delta_1     \|     u (t')    \|^2_{H^{1-s(t')}}  +
    \| v (t' )  \|^2_{H^{-s(t')}} \big) dt'
\end{equation}
\begin{equation}
\label{LH3}
+    \lambda    \int_0^t   e^{ -2 \gamma t'}
\big(  \delta_0 \delta_1     \|     u(t')    \|^2_{H^{1- s(t')+\ml} }  +
    \| v (t' )  \|^2_{H^{- s(t') +\ml} } \big) dt' .
\end{equation}

   On the right hand side, we find the initial data
\begin{equation}
\label{RH1}
C  A^2_{L^\infty}   \| u (0) \|^2_{H^{1- s(0)}}     +
      \| v(0) \|^2_{H^{-s(0)}}   ,
\end{equation}
the contribution of $f $
\begin{equation}
\label{RH2}
2 \re \int_0^t \la f  (t') ,   Q_\gamma v (t') \ra dt',
\end{equation}
an estimated contribution of $g$
\begin{equation}
\label{RH5}
   K  \int_0^t  e^{-2 \gamma t'}  \| g (t') \| _{1 - s(t) -  \ml }  \|
u (t' ) \|_{1 - s(t') + \ml }  dt'  ,
\end{equation}
and two types of   ``remainders'':
\begin{equation}
\label{RH3}
         K_0
      A_{LL} A_{L^\infty}  \int_0^t
       e^{-2 \gamma t'}    \big(   \| u  (t' ) \|^{2 }_{H^{1- s(t')+\ml
}}   +  \frac{1}{\delta_0^2}
   \| v_\gamma (t' )  \|^{2 }_{H^{ -s(t')+ \ml }} \big) dt'
\end{equation}
and
\begin{equation}
\label{RH4}
         K  \int_0^ t      e^{-2 \gamma t'}      \big(   \| u_\gamma(t'
) \|^{2 }_{H^{1- s(t')}}   +
   \| v(t' )  \|^{2 }_{H^{-s(t')}} \big)   dt' .
\end{equation}
If
\begin{equation}
\label{443}
\lambda \ge 2   K_0
      \frac{ A_{LL} A_{L^\infty} } {  \delta_0 \delta_1}   \quad \mathrm{and}
\quad  \lambda  \ge 2   K_0
      \frac{A_{LL} A_{L^\infty} }{\delta_0^2}
\end{equation}
the term in \eqref{RH3}  can be absorbed by \eqref{LH3}.
Note that this choice of $\lambda$ is precisely the choice announced
in \eqref{lda0},
with a new function $K_0$ of $ A_{L^\infty} / \delta_0$.
Finally, if $\gamma$ is large enough, the term \eqref{RH4} is
absorbed by \eqref{LH2},
    finishing  the proof of the main estimate \eqref{mainest22}.
\end{proof}

%%%%%%%%%%%%%%%%%%%%%%%%%%

\subsection{Proof of Theorem~\ref{thm42}}

 From now on, we assume that $\lambda \ge \lambda_0 $ and
 $\gamma \ge 
\gamma_0$ are fixed, so that the estimate
\eqref{mainest22} holds.
Consider $u$, $f$, $u_0$ and $u_1$ satisfying the equation \eqref{Cauchyb} and
the smoothness assumptions \eqref{assu1}, \eqref{assu2}, \eqref{assu3}.
Consider $v = X u + c_0 u$, which by Lemma~\ref{lem43} satisfies
\begin{equation}
\label{assu4}
v \in \cH_{- \theta + \ml}, \quad \D_t v \in L^1([0, T]; H^{-1- \theta_1}),
\quad  v_{| t = 0 } = v_0 \in H^{- \theta},
\end{equation}
with $v_0 = a_0 {}_{| t = 0}  u_1  + \sum  a_j{}_{| t =0} \D_{ x_j} u_0  +
c_0{}_{| t = 0} u_0 $.
In particular, $(u, v, f )$ and $g = 0$ satisfy \eqref{syste1}.

We mollify $u$ and $v $ and introduce, for $\eps > 0$,
\begin{equation}
u_\eps = J_\eps  u, \quad  v_{\eps} = J_ \eps   v  \quad \mathrm{with} \quad
J \eps = (1 - \eps \Delta_x)^{-1}.
\end{equation}
For all $\eps  > 0$, \eqref{assu1} and \eqref{assu4} imply that
\begin{eqnarray*}
u_\eps \in L^2([0, T], H^2), &\quad  \D_t u_\eps \in L^2([0, T], H^1),
\\
v_\eps \in L^2([0, T], H^1), &\quad  \D_t v_\eps \in L^1([0, T], L^2),
\end{eqnarray*}
(see \eqref{subassu}). Moreover, using the spatial Fourier transform,
one immediately
sees that $u_\eps $ converges to $u$ in  $\cH_{1- \theta, \lambda} (T)$
and  $v_\eps $ converges to $v$ in  $\cH_{- \theta, \lambda} (T)$.

Define
\begin{eqnarray*}
&& f_\eps   = Y^* v_\eps  + \tilde b_0 v_\eps  - \tilde L_2 u_\eps  +
   \tilde L_1 u_\eps  + \tilde d u_\eps,
   \\
&& g_\eps   = Y u_\eps  + \tilde c_0 u_\eps  - v_\eps / a_0  .
   \end{eqnarray*}

   \begin{lem}
   \label{lem46}
   Assumptions $\eqref{assu1}$ and $\eqref{assu3}$ imply that
   $f_\eps = f_{1, \eps} + f_{2, \eps}$ with
   $f_{1, \eps}  \to f_1$ in  $\cL_{- \theta, \lambda}(T)$ and
   $f_{2, \eps}  \to f_2$ in  $\cH_{- \theta - \ml , \lambda}(T)$.
   Moreover, $g_{ \eps}  \to 0$ in  $\cH_{1- \theta- \ml, \lambda}(T)$.
   \end{lem}

Taking this lemma for granted, we finish the proof of Theorem~\ref{thm42}.
We use the estimate \eqref{mainest22} for $(u_\eps, v_\eps)$,  together with
the estimates
$$
\begin{aligned}
| \la f_\eps (t) , Q_\gamma^2 v (t) \ra | \le
    C e^{- 2 \gamma t} \big( &   \| f_{1, \eps}(t)  \|_{H^{-s(t)}}
\|v_ \eps(t)  \|_{H^{-s(t)}} \\
   & +   \| f_{2, \eps}(t)  \|_{H^{-s(t)- \ml }}  \|v_ \eps(t)
\|_{H^{-s(t)+ \ml }} \big)
\end{aligned}
$$
and
$$
\begin{aligned}
\big| \int_0^t  \la f_\eps   , Q_\gamma^2 v   \ra  dt' \big| \le
    C\Big(  \int_0^t e^{-   \gamma t' } \big(    \| f_{1, \eps}(t')
\|_{H^{-s(t')}}  dt' \Big)
     \sup_{ 0 \le t' \le t} e^{- \gamma t'}   \|v_ \eps(t')  \|_{H^{-s(t')}} \\
    +  C
   \Big(  e^{- 2 \gamma t'}   \| f_{2, \eps}(t')  \|^2_{H^{-s(t')- \ml
}} dt'\Big)^\mez
    \Big(  e^{- 2 \gamma t'}   \|v_ \eps(t')  \|^2_{H^{-s(t')+ \ml }}
\big) dt'\Big)^\mez.
\end{aligned}
$$
This implies that there is a $K$ such that for all $\eps > 0$, one has
\begin{equation}
\label{mainest222}
\begin{aligned}
    \sup_{ 0 \le t' \le t} &  \| u_\eps  (t')  \|^2_{ H^{1 - s( t')}} +
    \sup_{ 0 \le t' \le t} \| v_\eps  (t')  \|^2_{ H^{-s(t') }}  \\
     +\int_0^t  \Big( &  \| u_\eps  (t' ) \|^2 _{  H^{1 - s( t') +
\frac{1}{2} log }} +
    \| v_\eps    (t' ) \|^2 _{  H^{ - s( t') + \frac{1}{2} log }} \Big) dt'
    \\
    \\
     \le    K  \Big\{ &  \| u _\eps(0)   \|^2_{ H^{1 - s(0) }} +
       \|  v_\eps(0)   \|^2_{ H^{-s(0) }}  + \int_0^t
\|g_\eps(t')\|^2_{H^{1- s(t) - \ml} } dt'
        \\
   &  +    \big( \int_0^t   \| f_{1, \eps}  (t' ) \|  _{  H^{- s( t')
}}  dt' \Big)^2  +
      \int_0^t   \| f_{2, \eps}  (t' ) \|^2 _{  H^{ -s( t')  -
\frac{1}{2} log }}  dt' \Big\} .
    \end{aligned}
\end{equation}
In addition, there are similar estimates for the differences
$(u_\eps - u_{\eps'}, v_\eps- v_{\eps'})$.
   Since $u_\eps (0) = J_\eps u_0 $ and $v_\eps (0) = J_\eps v_0 $ converge to
   $u_0$ and $v_0$ in $H^{1- s(0)}$ and  $H^{- s(0)}$, respectively,
   the estimate implies that
   $u_\eps $ is a Cauchy sequence in $\cH_{1- \theta, \lambda} (T)$ and
   in $C^0 ([0, t]; H^{1- s(t)})$ for all $t\in [0, T]$. Therefore, the limit
   $u$ in  $\cH_{1- \theta, \lambda} (T)$ also belongs to $\cC_{1-
\theta, \lambda} (T)$.
   Similarly, $v_\eps $ is a Cauchy sequence in $\cH_{- \theta,
\lambda} (T)$ and
   in $C^0 ([0, t]; H^{- s(t)})$ for all $t\in [0, T]$ and   $ v \in
\cC_{1- \theta, \lambda} (T)$.
In addition, we can pass to the limit in \eqref{mainest222} proving that
\begin{equation}
\label{mainest2222}
\begin{aligned}
    \sup_{ 0 \le t' \le t} &  \| u   (t')  \|^2_{ H^{1 - s( t')}} +
    \sup_{ 0 \le t' \le t} \| v   (t')  \|^2_{ H^{-s(t') }}  \\
     +\int_0^t  \Big( &  \| u   (t' ) \|^2 _{  H^{1 - s( t') +
\frac{1}{2} log }} +
    \| v   (t' ) \|^2 _{  H^{ - s( t') + \frac{1}{2} log }} \Big) dt'
    \\
    \\
     \le    K  \Big\{ &  \| u _0   \|^2_{ H^{1 - s(0) }} +
       \|  v_0   \|^2_{ H^{-s(0) }}
        \\
   &  +    \big( \int_0^t   \| f_{1}  (t' ) \|  _{  H^{- s( t')   }}
dt' \Big)^2  +
      \int_0^t   \| f_{2}  (t' ) \|^2 _{  H^{ -s( t')  -  \frac{1}{2}
log }}  dt' \Big\} .
    \end{aligned}
\end{equation}
Using the equation $Y u + \tilde c_0 u = v/ a_0$ and the estimate
\eqref{estv2} of Lemma~\ref{lem43} to
bound the time derivative $\D_t u$, we see that
$\D_t u \in \cC_{ - \theta, \lambda}(T) $ and that the energy estimate
   \eqref{mainest}  is satisfied.

Therefore, it remains only to prove the lemma.

\begin{proof}[Proof of Lemma~$\ref{lem46}$]
   By assumption \eqref{assu3},
$  f = f_1 + f_2 $ and $    J_\eps f_1  \to f_1 $  in  $ \cL_{-
\theta, \lambda}(T)$
and
$ J_\eps f_{2}  \to f_2 $   in   $ \cH_{- \theta - \frac{1}{2} log ,
\lambda}(T)$. Therefore, it
is sufficient to prove that the commutators
\begin{eqnarray*}
&&  [  Y^*, J_\eps]  v   ,  \quad    [ \tilde L_2 , J_\eps]  u
\\
&&  [ \tilde b_0 , J_\eps ] v,    \quad
[  \tilde L_1 , J_\eps ] u , \quad   [  \tilde d, J_\eps]  u ,
   \end{eqnarray*}
   converge to $0$ in $\cH_{- \theta, \lambda}(T)$
   and that the commutators
   $$
    [Y,   J_\eps] u  , \quad    [  \tilde c_0,  J_\eps] u, \quad    [ 1
/ a_0, J_\eps] v
   $$
   converge to $0$ in $\cH_{1- \theta, \lambda}(T)$.
   We note that $J_\eps$ commutes with $\D_t$ in $Y^*$ and
   $Y$. Thanks to \eqref{assu1} \eqref{assu4}
   and to the conservative form of $Y^*$ and $\tilde L*_2$, we see that
there are   four types
   of commutators  to consider :
\begin{equation}
\label{comm1}
\begin{aligned}
& [ a, J_\eps]  w  \to 0  \ \  \mathrm{in}\   \cH_{ 1 - \theta,
\lambda}(T), \quad   \mathrm{when } \\ &  a \in L^\infty \cap LL ([0,
T]\times \RR^d), \ \  w \in \cH_{- \theta, \lambda}(T),
\end{aligned}
\end{equation}
\begin{equation}
\label{comm2}
\begin{aligned}
& [ b, J_\eps]  w  \to 0  \ \  \mathrm{in}\   \cH_{  - \theta,
\lambda}(T), \quad   \mathrm{when } \\ &  b \in C^\alpha  ([0,
T]\times \RR^d), \ \  w \in \cH_{- \theta, \lambda}(T),
\end{aligned}
\end{equation}
    \begin{equation}
\label{comm3}
\begin{aligned}
& [ c, J_\eps]  w  \to 0  \ \  \mathrm{in}\   \cH_{ 1 - \theta,
\lambda}(T), \quad   \mathrm{when } \\ &  c \in C^\alpha  ([0,
T]\times \RR^d), \ \  w \in \cH_{1- \theta, \lambda}(T),
\end{aligned}
\end{equation}
      \begin{equation}
\label{comm4}
\begin{aligned}
& [ d, J_\eps]  w  \to 0  \ \  \mathrm{in}\   \cH_{  - \theta,
\lambda}(T), \quad   \mathrm{when } \\ &  d \in L^\infty ([0,
T]\times \RR^d), \ \  w \in \cH_{1- \theta, \lambda}(T).
\end{aligned}
\end{equation}

The first commutators $[ a, J_\eps] = [ T_a , J_\eps] +  R_a, J_\eps
- J_\eps R_a  $
are  uniformly bounded from
$ \cH_{- \theta, \lambda}(T)$ to $ \cH_{1- \theta, \lambda}(T)$: this
is true for the first
term by Proposition~\ref{propB3b}, since the $J_\eps$ form a bounded family of
operators of degree $0$; for the last two terms, this
follows from Proposition~\ref{proppara}.
Moreover,  $[a, J_\eps] w \to 0$ in $L^2([0, T]; H^{\sigma})$ for all
$\sigma < 1$, and thus also in
$\cH_{1- \theta, \lambda}$,
when $w$ is smooth and $a \in L\infty \cap LL$. By density, this
implies \eqref{comm1}.

For the commutators \eqref{comm2}, we note that they
are  uniformly bounded from
$ \cH_{- \theta, \lambda}(T)$ to $ \cH_{- \theta, \lambda}(T)$.
This is true for
both terms $ b  J_\eps  $   and $J_\eps b$ since $s(t) $ remains in a compact
subset of $[0, \alpha[$. Because $[ b, J_\eps]w $   converges to zero in
$L^2([0, T]; H^{\sigma})$ for all $\sigma < \alpha $,
when $w$ is smooth and $b \in L\infty \cap LL$, the convergence in
\eqref{comm2}
follows. The proof for \eqref{comm3} is similar.

Finally, we note that  $ [ d, J_\eps]  w  \to 0 $ in
$L^2([0, T] \times \RR^d) $, hence in $  \cH_{  - \theta, \lambda}(T)$
   when  $  d \in L^\infty ([0, T]\times \RR^d)$ and
   $ w \in L^2([0, T] \times \RR^d) $, thus in particular when $w \in
\cH_{1- \theta, \lambda}(T)$.
\end{proof}
%%%%%%%%%%%%%%%%%%%%%%%%

\subsection{Existence and uniqueness}

\begin{proof}[Proof of Theorem~$\ref{Uniq}$]
{\ }

Assume that $u \in H^{s} (] 0, T[ \times \RR^n ) $ with  $ s \in
]1-\alpha, \alpha [ $, $T \le T_0$,  and satisfies
\begin{equation}
\label{Cauchy0}
L u = 0, \quad  u_{|t = 0} = 0 , \quad  X u_{|t = 0} = 0.
\end{equation}
We want to  prove that $ u = 0 $.

Fix $ \theta < \theta_1$ in $ ]1-\alpha, \alpha [ $ with $1 - \theta < s$. Let
$ \lambda $ and $T' $  be the parameter and time associated to them
by  Theorem~\ref{thm42}.
Note that they depend  only on  $ \theta$,  $\theta_1$,
   the norms $A_{L^\infty} $ and $A_{LL}$ in \eqref{defA} and the
constants of hyperbolicity
$\delta_0$ and $\delta_1$ in \eqref{hyp}.

 From Lemma~\ref{lem42}, we know that
$ u \in L^2([0, T]; H^s(\RR^n))$ and
 $ \D_t u \in L^2([0, T]; H^{s-1}(\RR^n))$
and therefore, on $[0, T']\times \RR^n$,
$ u \in \cH_{1 - \theta + \ml, \lambda } $  and $ \D_t u \in \cH_{-
\theta + \ml, \lambda} $
since $s > 1 - \theta - \lambda t$.
By Theorem~\ref{thm42}, $u$ satisfies the energy estimate
\eqref{mainest} on $[0, T']$, and since the right hand side vanishes,
$u = 0$ for $ t < T'$.
By a  finite number of iterations, $u$ vanishes for $t < T$.
\end{proof}

\bigbreak

\begin{proof}[Proof of Theorem~$\ref{mainthm}$]{\ }

    On $[0, T_0] \times \RR^d$,  the coeffficients of $L_2 $ can be
approximated in  $L^\infty$ and  $C^{\alpha'}$  for all $\alpha' < 1$
by
$C^\infty$ functions which are uniformly bounded in $L^\infty $ and
in $LL $, in such a way that the hyperbolicity condition
\eqref{hyp}  remains satisfied.
Similarly, the coeffficients of $L_1 $ can be approximated in
$L^\infty$ and  $C^{\alpha'}$  for all $\alpha' < \alpha$ by
smooth functions which are uniformly bounded in $ C^\alpha$. Further,  the
coefficient $c$ can be approximated
in $L^2_{loc} $ by functions uniformly bounded in $L^\infty$.
This 
defines operators  $L^\eps$ with $C^ \infty $ coefficients which
satisfy \eqref{hyp},
\eqref{defA} and \eqref{defB} uniformly in $\eps$ and converge to the
coefficients of $L$ in the sense
described above.

We fix the parameter $\lambda \ge  \lambda_0$, where $\lambda_0$ is given by
Theorem~\ref{thm42}. Recall that $T$ is then given by \eqref{fixT}.
Consider Cauchy data $u_0 \in H^{ 1- \theta}$ and $u_1 \in H^{ -
\theta}$ and  a
source term  $f = f_1 + f_2 $
with $f_1 \in \cL_{- \theta, \lambda} (T) $ and $ f_2 \in \cH_{-
\theta-\ml, \lambda}$.  We can approximate these data in the
corresponding  spaces
by $C^\infty $  functions $u_0^\eps$, $u_1^\eps$ ,  $f_1^\eps$ and
$f_2^\eps$, compactly supported in $x$.
   The strictly hyperbolic problems with smooth coefficients and smooth data
\begin{equation}
\label{regeq}
L^ \eps u^ \eps = f_1^\eps+ f_2^\eps , \quad u^\eps{}_{|t = 0} = u_0^\eps ,
\quad X^\eps u^\eps{}_{|t = 0} = u_1^\eps
\end{equation}
have a unique smooth solution $u^\eps$, compactly supported in $x$.

By Theorem~\ref{thm42}, the energy estimate \eqref{mainest} is satisfied
with a constant $K$ independent of $\eps$. Therefore the family
$\{ u^\eps\}$ is bounded in $ \cH_{1- \theta+ \ml, \lambda}$, thus in
$L^2([0, T], H^{1-\theta_1})$  and
the families $\{ \D_t u^\eps\}$  and  $\{ X^\eps u^\eps\}$ are bounded in
$ \cH_{- \theta+ \ml, \lambda}$, hence   in
$L^2([0, T], H^{-\theta_1})$.
Therefore, extracting a subsequence if necessary,  $u^\eps$ converges
to a limit $u$,  weakly  in
$L^2([0, T], H^{1-\theta_1}) $  and in
$H^1([0, T], H^{-\theta_1}) $.  Moreover,
$ u \in  \cH_{1- \theta+ \ml, \lambda}$ and
    $ \D_t u\in  \cH_{- \theta+ \ml, \lambda}$.
There is  no difficulty in passing to the limit in the equation in
the sense of distributions:
all the products are well defined and involve  one strong and one
weak convergence.
Thus $L u = f$.

The weak convergence in  $L^2([0, T], H^{1-\theta_1}) \cap H^1([0,
T], H^{-\theta_1}) $
implies the strong convergence in $ C^0([0, T] ; H^{ - \theta_1}_{loc} ) $
and therefore the convergence of   $u^\eps_{|t = 0} $ to  $u_{|t =
0}$ in $H^{ - \theta_1}_{loc} $.
Therefore,  $u_{|t = 0} = u_0$.

Using the equation as in Lemma~\ref{lem42}, we prove that  the family
$ v^\eps = X^\eps u^\eps + c_0^\eps u^ \eps$, which converges weakly to
$v = X u + c_0 u$, is bounded in
$L^2([0, T], H^{-\theta_1}) \cap H^1([0, T], H^{-1 -\theta_1}) $. Thus
$v^\eps_{|t = 0} $  converges to  $v_{|t = 0}$ in $H^{ - \theta_1}_{loc} $.
Hence  $ v_{|t = 0} = u_1 + c_0{}_{|t = 0}  u_0 $ implying that
$ X u_{|t = 0 } = u_1$.

By Theorem~\ref{thm42} the solution  $u$  also belong to
$\cC_{1-\theta, \lambda}$ with $ \D_t u\in   \cC_{1-\theta, \lambda}$
and satisfies the energy estimate \eqref{mainest}.
\end{proof}

%%%%%%%%%%%%%%%%%%%%%%%%%%
%%%%%%%%%%%%%%%%%%%%%%%%%%
%%%   SECTION 5

\section{Local results }

We consider the equation \eqref{maineq} together with an initial
hypersurface $\Sigma$
satisfying Assumption~\ref{assl}.
Everything being local, and the assumptions being invariant under
smooth changes  of
coordinates, we may assume that  we are working in  coordinates
$ y = (t, x)$ such that   $ \underline y = (0, 0)$ and $\Sigma = \{ t = 0 \}$.
The operator has the form  \eqref{eq1} with coefficients which are defined
on a neigborhood  $ \Omega$ of the origin.

Lemma~\ref{traces}  is a local version of
Lemmas~\ref{lem41} and \ref{lem42}. The proof is identical, using
local multiplicative
properties and local versions of the spaces $H^{s, s'}$.

\subsection{Local existence}
Choose  $ \Phi $, a smooth map from
$\RR^{1+n}$ to $\Omega$,  with $ \Phi (y) = y $ on a  smaller
neighborhood $\Omega_1$ and
$ \Phi ( y) = 0 $  for $ y$ large enough.
Changing the coefficients acoording to the rule $ a^\sharp (y) =  a( \Phi(y))$
we obtain an operator $L^ \sharp $ which coincides with $L$ on
$\Omega_1$, satisfies
the regularity conditions \eqref{cf1} to \eqref{cf3}, and the
hyperbolicity conditions \eqref{hyp}
globally on $ \RR^{1+ n}$.

Fix  $ s > 1 - \alpha$. Without loss of generality for the statememt
of Theorem~\ref{Th1}, we can assume that $ s < \alpha$. We are going
to apply Theorem~\ref{mainthm}
to the operator $L^\sharp$ with $ \theta = 1 - s \in ] 1 - \alpha,
\alpha [$.  Choosing
$ \theta_1 \in ]\theta, \alpha[ $, this theorem provides us with
$\lambda $ and  $T =  ( \theta_1- \theta)/ \lambda$.
We fix  $ \Omega'  = \Omega_1 \cap \{ | t | < T \}$.

Suppose that  $u_0$ and $u_1$ are Cauchy data in $ H^s (\omega)$ 
and 
$ H^{s-1} (\omega)$  respectively,
on a neighborhood $\omega$ of $0$ in $\RR^n$.
    There are restrictions to $\omega$ of functions
   $u^\sharp_0 \in H^s (\RR^n)$ and $u^\sharp_1 \in  H^{s-1} (\RR^n)$
respectively.
Suppose that  $ f \in L^2 ( \Omega' \cap \{ t > 0\} ) $. We extend
it, for instance by $0$,
to $ f^\sharp \in L^2 ( [ 0, T] \times \RR^n)$.
By Theorem~\ref{mainthm}, the Cauchy problem
\begin{equation}
\label{CPxx}
L^\sharp  u^\sharp  =  f^\sharp  , \quad  u^\sharp _{|  t =0 } =  u^\sharp_0,
\quad  (X^\sharp  u^\sharp )_{| t = 0 } = u^\sharp _1
\end{equation}
has a solution  $ u^\sharp $ on $[0, T] \times \RR^n$, which belongs
in particular
to $ L^2 ( [0, T]; H^ {s_1}) $  with   $ s_1 = 1 - \theta_1$  and such that
$ \D_t u \in L^2 ( [0, T]; H^ {s_1 -1}) $. In particular,
$u^\sharp \in H^ { s_1} ( [0, T] \times \RR^n)$ and by restriction
to $\Omega '$ defines a solution of \eqref{CP}.

\subsection{Local uniqueness}

To prove Theorem~\ref{Th2},
we first reduce  the problem  to proving a theorem of propagation of
zero accross the surface $ \{ t = 0 \}$.

\begin{lem}
\label{lem51}
   Suppose  that  $ s > 1- \alpha$ and $u \in H^s ( \Omega \cap \{ t >
0 \} ) $ satisfies
   \begin{equation}
\label{CP0b}
L u = 0   , \quad  u_{| t = 0 } =  0, \quad  X  u_{| t = 0 } = 0.
\end{equation}
Then the extension $  u_e $  of $u$ by  $0$ for $ t < 0 $
satisfies
\begin{equation}
\label{exteq}
u_e \in H^s \quad \mathrm{and} \quad  L   u_e  = 0
\end{equation}
on a neighborhood  $\Omega_1$ of $0$.
\end{lem}

\begin{proof}
If the coefficients were smooth, this would be immediate. We check that
we have enough smoothness
to extend the result to our case.

We can assume that  $ \Omega = ] - T, T [ \times \omega $.   From
Lemma~\ref{lem42}  (localized in space)  we know
that  $ u \in L^2( [0, T] ; H^ s_{loc}(\omega))$, thus  its extension
   $ u_e \in L^2( [-T, T] ; H^ s_{loc}(\omega))$. Moreover,
   $ \D_t u \in L^2( [0, T] ; H^ {s-1}_{loc}(\omega))$  and by assumption
$ u_{| t = 0 } =  0$. Therefore,
$\D_t u_e $ is the extension of $\D_t u $ by $0$ and thus belongs to
$ L^2( [-T, T] ; H^{ s-1} _{loc}(\omega))$.
In particular,  $u_e \in H^s_{loc} ( ] - T, T [ \times \omega)$.

Let $ v = Xu + c_0 u  \in L^2( [0, T] ; H^ {s-1}_{loc}(\omega))$ and
let $ v_e \in L^2( [-T, T] ; H^ {s-1}_{loc}(\omega))$  denote  its
extension by $0$.
The first step  implies that $X u_e$  is the extension of  $ Xu $ and therefore
$ v_e  = X u_e + c_0 u_e$.
Write the equation as
\begin{equation}\label{53}
\D_t v = P (u, v)
\end{equation}
where $P$ involves only spatial derivatives (see \eqref{syste2}). Morever,
we have seen in the proof of Lemma~\ref{lem42} that
$ P(u, v) \in   L^2( [0, T] ; H^ {s-2}_{loc}(\omega))$.
Since by assumption the trace of  $v$ vanishes, this implies that
$ \D_t v_e $ is the extension by $0$ of $\D_t v$, thus the extension of
$P(u, v)$,  that is $ P( u_e, v_e)$. Since $ v_e  = X u_e +c_0 u_e$,
this  means that $u_e$  satisfies the equation on
$ \Omega = ] - T, T [ \times \omega$.
\end{proof}

We now finish the proof of Theorem~\ref{Th2}.
We  suppose that  $u \in H^s ( \Omega \cap \{ t > 0 \} ) $ satisfies
\eqref{CP0b}, with  $ s > 1- \alpha$ and
we denote by $u_e$ its  extension   by  $0$ for $ t < 0 $.
We use the classical convexification method, and consider the change
of variables
\begin{equation}
\label{54}
(t, x) \mapsto  (\tilde t , \tilde x)     \quad   \tilde t =    t + |
x|^2 , \  \tilde x = x ,
\end{equation}
which maps  the past $ \{ t < 0 \} $ to  $ \{ \tilde t <    | \tilde
x |^ 2 \} $.
Thus there is $ T_0 > 0$ such that the  function $\tilde u$ deduced from $u_e$
is defined for $ \tilde t < T_0 $ and vanishes for  $ \tilde t <  |
\tilde x |^ 2 $.
Moreover, decreasing $T_0$ if necessary, the operator
$ \tilde L$ deduced from  $L$  is defined
on a neighborhood $\tilde  \Omega$  of the origin which contains  the
closed  lens
$ \overline D = \{ | \tilde x|^2  \le t \le T_0 \} $
and $\tilde L \tilde u = 0 $ on $ \tilde \Omega \cap \{ t < T_0 \} $.
Now we extend the coefficients of $\tilde L$, as above, and obtain a
new  operator
$ L^\sharp $, defined on $ \RR^{1+ n}$, satisfying the assumptions of
section 2, and
equal to $\tilde L$ on a neighborhood of $ \overline D $.
Therefore,  on $ ]- \infty, T_0[ \times \RR^n $
\begin{equation}
\label{PL}
L^\sharp \tilde u = 0, \quad  \tilde u \in H^ s   , \quad  \tilde
u_{| \{ \tilde t < | \tilde x|^2 \} }
= 0.
\end{equation}
Since $ \tilde u$ vanishes in the past, the traces $\tilde u_{| t = -
\eps} $  and
$X^\sharp \tilde u_{| t = - \eps} $  vanish for all $ \eps > 0$.
Therefore, Theorem~\ref{Uniq} applied to the Cauchy problem for
$L^\sharp$  with initial time $- \eps$ implies that $\tilde u = 0 $
for all $(\tilde t, \tilde x) $  such that  $ \tilde t < T_0$.
Hence  $u = 0$  on a neighborhood of the origin.

   %%%%%%%%%%%%%%%%%%%%

%%%%%%%%%%%%%%%%%%%%%%%%%%

\end{document}